\documentclass[a4paper,11pt,reqno]{article}
\usepackage{amsmath, amsfonts, amssymb, amsthm}
\usepackage{mathrsfs}
\usepackage{amsmath}
\usepackage{amssymb,color}
\usepackage{graphicx}
\usepackage{amssymb}

\usepackage[colorlinks=true]{hyperref}
\hypersetup{urlcolor=blue, citecolor=blue}

\par

\def \R {{\mathbb{R}}}
\numberwithin{equation}{section}

\begin{document}

\title{Two-component system modelling shallow-water waves with constant vorticity under the Camassa-Holm scaling}

\author{Leyi $\mbox{Zhang}$\footnote{E-mail:
TS22080025A31@cumt.edu.cn}\quad
and \quad Xingxing $\mbox{Liu}$\footnote{Corresponding author. E-mail: liuxxmaths@cumt.edu.cn}\\
$\mbox{School}$ of mathematics, China University of
Mining and Technology,\\Xuzhou, Jiangsu 221116, China}

\date{}
\maketitle

\begin{abstract}
This paper is concerned with the derivation of a two-component system modelling shallow-water waves with constant vorticity under the Camassa-Holm scaling
from our newly established Green-Naghdi equations with a linear shear. It is worth pointing out that the $\rho$ component in this new system is quite different from the previous two-component system
due to the effects of both vorticity and larger amplitude. We then establish the local well-posedness of this new system in Besov spaces, and  present a blow-up criterion.
We finally give a sufficient condition for global strong solutions to the system in some special case.\\


\end{abstract}

\section{Introduction}
It is well-known that compared to the KdV equation, the great interest in the Camassa-Holm (CH) equation lies in the fact that it could exhibit both the phenomena of wave breaking and soliton interaction \cite{Brandolese1,C-H,Constantin3,Constantin,Constantin-E,Constantin-L,Constantin-S1,Constantin-S2,Johnson2}. In the past few years, various multi-component integrable system of coupled equations have been proposed to generalize the CH equation.
The one that has caught the most attention is the following two-component CH system \cite{CLZ,Constantin-I,Ivanov,Olver-R}
\begin{equation}\label{2CH}
\left\{
 \begin{aligned}
&{m_t}-Au_x+2u_xm+um_x+\rho\rho_x=0,\\
&{\rho _t}+ (\rho u)_x=  0,
\end{aligned}
\right.
\end{equation}
where $m=u-u_{xx}$. The variable $u(t,x)$ represents the horizontal velocity of the fluid, and $\rho(t,x)$ is in connection with the free surface deviation from equilibrium
with the boundary assumptions $u\rightarrow 0$ and $\rho\rightarrow1$ as $|x|\rightarrow \infty$. The parameter $A\geq 0$ is a constant which describes the vorticity of the underlying flow.
The system (\ref{2CH}) without vorticity, i.e. $A=0$, was firstly derived by Constantin and Ivanov \cite{Constantin-I} in the context of shallow water theory. To incorporate an underlying vorticity in the flow, Ivanov \cite{Ivanov}
gave a justification of the derivation of the system (\ref{2CH}). Then a quite extensive study of the qualitative properties for system (\ref{2CH}) was carried out in a lot of works, e.g.
\cite{CLQ,Escher-L-Y,G-Y,G-L1,G-L,L,ZY}.

Recently, Hu and Liu \cite{H-L} established the following two-component CH system modelling shallow-water waves propagating in the equatorial ocean regions
with the Coriolis effect due to the Earth's rotation
\begin{equation}\label{coriolis}
\left\{
 \begin{aligned}
&{\rho _t} + {(\rho u)_x} + \Omega\rho {({u^2})_x} + 8{\tilde{\beta} _1}\rho {({u^3})_x} =  0,\vspace{1ex}\\
&{m_t} + \sigma (2m{u_x} + u{m_x}) + 3(1 - \sigma )u{u_x} + \frac{1}{2}{({\rho ^2})_x}+2\Omega(\rho ^2 u)_x-8\Omega(\rho u)_x
\vspace{1ex}\\
 &\qquad\qquad\qquad\qquad\qquad\qquad+4\Omega u_x + 24{\tilde{\beta} _1}{({u^2}\rho (\rho  - 1))_x} + 4{\tilde{\beta} _2}{({u^3})_x} = 0,
\end{aligned}
\right.
\end{equation}
where $u(t,x)$ is the average horizontal fluid velocity, $m = u - {u_{xx}}$ is the momentum density, and $\rho(t,x)$ is also related to the free surface elevation from equilibrium.
The constants ${\tilde{\beta} _1}= \frac{{1-7{c^2}}}{{48c^2}},$ ${\tilde{\beta} _2}= \frac{{{c^4} + 5c^2-2}}{4c^3},$ where $c=\sqrt{1+\Omega^2}-\Omega$, with the parameter $\Omega$ which is the constant rotational frequency due to the Coriolis effect. The dimensionless parameter $\sigma$ shows a balance between nonlinear steepening and amplification in fluid convection due to stretching.
It is noted that the establishment of the system (\ref{coriolis}) is different from (\ref{2CH}) mainly based on the following two aspects. On the one hand,
the Coriolis force is not negligible when considering water wave motions which are presented on a huge range of spatial and temporal scales. The ocean dynamics near the Equator is quite different from that in non-equatorial regions since the change of sign of the Coriolis force across the Equator produces an effective waveguide, with the Equator acting alike to a natural boundary that facilitates azimuthal flow propagation \cite{Constantin-I1,Constantin-J2}. On the other hand, in the derivation of simple shallow-water models, we usually need two fundamental dimensionless parameters: the amplitude parameter $\varepsilon:=\frac{a}{h_{0}}$ and the shallowness parameter $\mu:=\frac{h_{0}^{2}}{\lambda^{2}},$ where $h_0$ is the mean depth of water, $a$ and $\lambda$ are the typical amplitude and wavelength of the waves, respectively.
In the shallow water regime ($\mu\ll 1$), with additional assumption that $\varepsilon\ll 1$, many interesting
water wave models of the corresponding specific regime can be derived from the governing
equation by making assumptions on the respective size of $\varepsilon,\mu$ \cite{Constantin-I,Constantin-J,Johnson1,Lannes}. For our situation here,
the system (\ref{2CH}) was derived from the Green-Naghdi \cite{G-N} in the Boussinesq scaling (weakly nonlinear regime) $\mu \ll 1$, $\varepsilon=O(\mu)$, whereas
the system (\ref{coriolis}) was obtained from the rotation-Green-Naghdi equations \cite{G-L-L} in the CH scaling (moderately nonlinear regime) $\mu \ll 1$, $\varepsilon=O(\sqrt{\mu})$.
Thus we can see that the larger the parameter $\varepsilon$, the stronger the nonlinear effects. Hu and Liu \cite{H-L} considered
the local well-posedness for (\ref{coriolis}) with initial data $(u_0,\rho_0)\in H^s\times H^{s-1}$ with $s>\frac{3}{2}$ and blow-up criteria. Moreover, they investigated the wave-breaking phenomena and existence of
global solution for (\ref{coriolis}) with $\Omega=0$ and $\sigma=1$. Dong \cite{Dong} gave some new wave breaking criteria and extended the earlier blow-up results of (\ref{coriolis}) with $\Omega=0$ and $\sigma\neq0$.

In the present paper, inspired by \cite{Constantin-I,H-L}, we derive a new two-component system modelling shallow-water waves with constant vorticity under the Camassa-Holm scaling in Section \ref{sec2}.
To this end, we start from the derivation of the Green-Naghdi equations with constant vorticity
from the governing equations for two-dimensional water waves in the presence of shear. Then introducing a new variable $\rho$ as in \cite{Constantin-I}, we formally obtain the following
two-component system with an underlying constant vorticity in the CH scaling
\begin{equation}\label{2CHvor}
\left\{
 \begin{aligned}
&{\rho _t} + \frac{1}{3}A{({\rho ^3})_x} + {(\rho u)_x} + \frac{1}{2}A\rho {({u^2})_x} + 8{\beta _1}\rho {({u^3})_x} =  0,\vspace{1ex}\\
&{m_t} + \sigma (2m{u_x} + u{m_x}) + 3(1 - \sigma )u{u_x} + \frac{1}{2}{({\rho ^2})_x}\vspace{1ex}\\
 &\qquad\qquad\qquad + 24{\beta _1}{({u^2}\rho (\rho  - 1))_x} + 4{\beta _2}{({u^3})_x} - A{u_{xxx}} = 0,
\end{aligned}
\right.
\end{equation}
where $u(t,x)$ is the average horizontal fluid velocity, $m = u - {u_{xx}}$ is the momentum density, and $\rho(t,x)$ is related to the free surface elevation from equilibrium.
The coefficients ${\beta _1}= - \frac{{{c^2} + 5}}{{48}},$ ${\beta _2}= \frac{{{c^3} + 3c}}{4},$ and $c=\frac{A+\sqrt{A^{2}+4}}{2}$, where $A$ represents the vorticity of the underlying flow.
The parameter $\sigma$ also shows a balance between nonlinear steepening and amplification in fluid convection due to stretching.
Note that the system $(\ref{2CHvor})$ with $A=0$ is exactly same as the one (\ref{coriolis}) with $\Omega=0.$ That is
\begin{equation}\label{two-component0}
\left\{
 \begin{aligned}
&{m_t} + \sigma (2m{u_x} + u{m_x}) + 3(1 - \sigma )u{u_x} + \frac{1}{2}{({\rho ^2})_x}\vspace{1ex}-3{({u^2}\rho (\rho  - 1))_x} + 4{({u^3})_x} = 0\vspace{1ex},\\
&{\rho _t}  + {(\rho u)_x} -\rho {({u^3})_x} =  0,
\end{aligned}
\right.
\end{equation}
where $m = u - {u_{xx}}.$ The system (\ref{two-component0}) can be regarded as a generalized two-component CH system under the Camassa-Holm scaling.

In the sequel, for notational convenience, we shall consider the following Cauchy problem of the system of (\ref{2CHvor}) with more general coefficients
\begin{equation}\label{two-component2}
\left\{
 \begin{aligned}
&{m_t} + \sigma (2m{u_x} + u{m_x}) + 3(1 - \sigma )u{u_x} + \frac{1}{2}{({\rho ^2})_x}+a_1(\rho^2u)_x+a_2(\rho u)_x\vspace{1ex}\\
 &\qquad\qquad\qquad \qquad\qquad+a_3u_x+a_4{u_{xxx}}+a_5{({u^2}\rho (\rho  - 1))_x}+a_6{({u^3})_x}= 0,\\
&{\rho _t} +{(\rho u)_x}+b_1 \rho^2 \rho_x+b_2 \rho u u_x+b_3\rho u^2 u_x =  0,\vspace{1ex}
\end{aligned}
\right.
\end{equation}
where $a_i,b_i$ are real constants. Obviously, the systems (\ref{coriolis}) and (\ref{2CHvor}) are the special cases of (\ref{two-component2}). Denoting $\zeta:= \rho -1$, thus $\zeta\rightarrow0$ as $|x|\rightarrow\infty$.
Then applying the translations $u(t,x)\rightarrow u(t,x+a_4t),$   $\zeta(t,x)\rightarrow u(t,x-b_1t)$ to the system (\ref{two-component2}), for $t>0,x\in \R$, we have
\begin{equation}\label{two-component}
\left\{
 \begin{aligned}
&u_t+\sigma uu_x=P(D)\big(\frac{3-\sigma}{2}u^2+\frac{\sigma}{2}(u_x)^2+\frac{1}{2}\zeta^2+\zeta+a_1(\zeta^2 u+2\zeta u+u)\\
&\qquad+a_2(\zeta u+u)+(a_3+a_4)u+a_5(u^2\zeta^2+u^2\zeta)+a_6u^3\big),\\
&\zeta_t+u\zeta_x+b_1\zeta^2\zeta_x+2b_1\zeta\zeta_x=-\zeta u_x-u_x-b_2\zeta uu_x-b_2uu_x-b_3\zeta u^2u_x-b_3u^2u_x,\\
&u(0,x)=u_0(x),\\
&\zeta(0,x)=\zeta_0(x)=\rho_0(x)-1,
\end{aligned}
\right.
\end{equation}
where $P(D)$ is the Fourier integral operator with the Fourier multiplier $-i\xi(1+\xi^2)^{-1}$. In Section \ref{sec3}, we show the local well-posedness of the Cauchy problem of the general system (\ref{two-component}) in nonhomogeneous Besov spaces. And thus it implies that (\ref{two-component}) is locally well-posed in
$H^s(\R)\times H^{s-1}(\R)$, $s>\frac{5}{2}.$ In comparison with (\ref{coriolis}),
the variable $\zeta$ in $(\ref{two-component})$ is "transported" along direction of $(u+b_1\zeta^2+2b_1\zeta)$ instead of $u$ in $(\ref{coriolis})$, which leads to
more regularity assumptions on initial data for the local well-posedness of $(\ref{two-component})$. Another difficulty is to obtain the estimate of the uniform bound of the approximate solutions, which results from
the higher-order nonlinearities appearing in $(\ref{two-component})$. Here we deal with it by using the mean-value theorem of integrals and choosing reasonably the
existence time of solution. In Section \ref{sec4}, applying commutator estimates in the nonhomogeneous Besov framework, we present a precise blow-up scenario for the general system (\ref{two-component}).
Here we can no longer control $\|\rho\|_{L^\infty}$ in terms of $\|u_x\|_{L^\infty}$ as done for the $\rho$ component in $(\ref{coriolis})$, since
$\zeta$ in $(\ref{two-component})$ is advected by $(u+b_1\zeta^2+2b_1\zeta)$. Finally, in Section \ref{sec5}, we give a sufficient condition for global
strong solutions to (\ref{two-component0}) in the special case of $\sigma=0$.

\section{Derivation of the model system}\label{sec2}
\newtheorem{theorem2}{Theorem}[section]
\newtheorem{lemma2}{Lemma}[section]
\newtheorem {remark2}{Remark}[section]
\newtheorem {definition2}{Definition}[section]
\newtheorem{corollary2}{Corollary}[section]
\par
In this section, we first briefly recall the governing equations in the presence of shear. Then
we establish the Green-Naghdi equations with a linear shear. Finally, in the Camassa-Holm scaling, we derive a two-component shallow-water wave system with an underlying constant vorticity
from the new obtained Green-Naghdi equations with shear.

\subsection{Governing equations in the presence of shear}
To start with, we consider the following the governing equations for the incompressible, inviscid and two-dimensional fluid motion (see, e.g.,\cite{Johnson1,Lannes}),
which consist Euler's equations and the equation of mass conservation
\begin{align}
\nonumber
\begin{cases}
u_{t}+uu_{x}+wu_{z}=-\frac{1}{\rho}p_{x}, &$in$ \ D_{t},\\
w_{t}+uw_{x}+ww_{z}=-\frac{1}{\rho}p_{z}, &$in$ \ D_{t},\\
u_x+w_z=0, &$in$ \ D_{t},
\end{cases}
\end{align}
together with the dynamic and kinematic boundary conditions
\begin{align}
\nonumber
\begin{cases}
p=\rho g\eta, &$on$ \ z=h_{0}+\eta(t,x),\\
w=\eta_{t}+u\eta_{x}, &$on$ \ z=h_{0}+\eta(t,x),\\
w=0,  &$on$ \ z=0,
\end{cases}
\end{align}
where $h_0$ is the mean depth, or the undisturbed depth of water, $\eta(t,x)$  measures the deviation from the average level, and the domain $D_{t}=\{(x,z):0<z<h_{0}+\eta(t,x)\}$.
$(u(t,x,z),w(t,x,z))$ is the two-dimensional velocity field, $g$ is the constant
Earth's gravity acceleration, and the variable $p$ measures the deviation from the hydrostatic pressure distribution.

Then one can perform the nondimensionalisation procedure on the above governing equations using the physical characteristic of the flow.
This process relies on the the amplitude parameter $\varepsilon=\frac{a}{h_{0}}$ and the shallowness parameter $\mu=\frac{h_{0}^{2}}{\lambda^{2}}.$ In order to incorporate waves in the presence of a shear flow, we consider an exact solution of the form $(u,w,p,\eta)\equiv (U(z),0,0,0)$, which
represents laminar flows with a flat free surface and with an arbitrary shear. Finally, making use of the transformation of scaling $u\rightarrow U(z)+\varepsilon u, w\rightarrow\varepsilon w,  p\rightarrow\varepsilon p,$
the governing equations in terms of dimensionless variables in the presence of a shear flow become (we refer to \cite{Escher-H-K-Y,Ivanov,Johnson4,L3} for the details)
\begin{align}\label{scaling and shear}
\begin{cases}
u_{t}+Uu_x+wU'+\varepsilon(uu_{x}+wu_{z})=-p_{x}, &$in$ \  0<z<1+\varepsilon\eta(t,x),\\
\mu(w_{t}+Uw_x+\varepsilon(uw_{x}+ww_{z}))=-p_{z}, &$in$ \ 0<z<1+\varepsilon\eta(t,x), \\
u_{x}+w_{z}=0, &$in$ \ 0<z<1+\varepsilon\eta(t,x),\\
p=\eta, &$on$ \ z=1+\varepsilon\eta(t,x),\\
w=\eta_{t}+(U+\varepsilon u)\eta_{x}, &$on$ \ z=1+\varepsilon\eta(t,x),\\
w=0, &$on$ \ z=0.
\end{cases}
\end{align}

For the nontrivial special case of a laminar shear flow, we taking $U(z)=Az$. The familiar Burns condition \cite{Burns} gives an expression for the speed $c$ of the travelling waves in linear approximation:
$\int_0^1\frac{dx}{(Az-c)^2}=1,$ which implies $c=\frac{1}{2}(A\pm\sqrt{A^{2}+4})$. Before the scaling the vorticity $\omega=U'+u_z-w_x.$ Using the nondimensionalised variables and the scalings of $u,w$ as before, and $\omega\rightarrow\ \sqrt{g/h_0}\omega$, we get $\omega=A+\varepsilon(u_z-\mu w_x)$. Note that the Burns condition arises as a local bifurcation condition for shallow water waves with constant vorticity \cite{Constantin0}. Thus, to find a solution with constant vorticity, we obtain
\begin{eqnarray}\label{vorticity}
u_z-\mu w_x=0,
\end{eqnarray}
which yields that the vorticity $\omega\equiv A.$ Therefore, we find that Eqs. (\ref{scaling and shear})-(\ref{vorticity}) are our final version of the governing equations under study.

\subsection{Derivation of the Green-Naghdi equations with a linear shear}
In this subsection, inspired by the derivation of the rotation-Green-Naghdi equations \cite{G-L-L},
we apply formal asymptotic methods to the full governing equations (\ref{scaling and shear})-(\ref{vorticity}) with $U(z)=Az$ to drive the Green-Naghdi equations with a linear shear in the shallow-water scaling ($\mu\ll1$), without any assumption on $\varepsilon$ (that is, $\varepsilon=O(1)$).

For the first equation of the Green-Naghdi equations with a linear shear, we denote $\bar{u}$ as the vertically averaged
horizontal component of the velocity,
\begin{equation}\label{ubar}
\bar{u}(t,x)=\frac{1}{h}\int_{0}^{h}u(t,x,z)\mathrm{~d}z,
\end{equation}
where $h=h(t, x)=1+\varepsilon \eta(t, x)$. Multiplying both sides of (\ref{ubar}) by $h$ and differentiating with respect to $x$, we obtain
$$(h \bar{u})_x=\int_0^h u_x \mathrm{~d} z+\varepsilon \eta_x u(t,x,h).$$
In view of the third equation in Eq. (\ref{scaling and shear}), substituting $u_x$ by $-w_{z}$ and taking account of the boundary conditions of $w$ on $z=0$ and $z=h$, we have
\begin{equation}\label{eta}
\eta_t+Ah\eta_x+(h\bar{u})_x=0,
\end{equation}
or
\begin{equation}\label{h}
h_t+Ahh_x+\varepsilon(h\bar{u})_x=0.
\end{equation}

For the second equation of the Green-Naghdi equations with a linear shear, we divide the process of derivation into two parts.
Assume that
\begin{equation}\label{u}
u(t,x,z) = {u_0}(t,x,z) + \mu {u_1}(t,x,z) + O({\mu ^2}),
\end{equation}
and
\begin{equation}\label{w1}
w(t,x,z) = {w_0}(t,x,z) + \mu {w_1}(t,x,z) + O({\mu ^2}).
\end{equation}
Firstly we get the expressions of $u$, $w$, $p$ in terms of $u_0$, $h$, $z$, and find the equation related to $u_0$ and $h$ only.
Combining (\ref{vorticity}) with (\ref{u})-(\ref{w1}), we obtain the following equalities in asymptotic order of $\mu^0$-level and $\mu^1$-level:
\begin{equation}\label{01level}
{\mu ^0}:\ u_{0,z} = 0\quad  \mbox{and} \quad {\mu ^1}:\ w_{0x}=u_{1z}.
\end{equation}
Thus the first equality of (\ref{01level}) implies that $u_0$ is independent of $z$, i.e., $u_0 = u_0(t,x)$. On the other hand,
it follows from the third equations in (\ref{scaling and shear}) and (\ref{u})-(\ref{w1}) that
\begin{equation}\label{01level0}
{\mu ^0}:\ {w_{0,z}}=-{u_{0,x}}\quad  \mbox{and} \quad {\mu ^1}:\ {w_{1,z}}=-{u_{1,x}}.
\end{equation}
Due to $w=0$ on $z=0$, we deduce from integrating the first equality of (\ref{01level0}) that
$w_0=- z{u_{0,x}}.$ Then integrating the second equality of (\ref{01level}) with respect to $z$ yields ${u_1} = -\frac{z^2}{2}{u_{0,xx}}.$
Moreover, it follows from the the second equality of (\ref{01level0}) that $w_1=\frac{{{z^3}}}{6}{u_{0,xxx}}.$ Therefore, we have
\begin{equation}\label{finalu}
u = {u_0} - \mu \frac{{{z^2}}}{2}{u_{0,xx}} + O({\mu ^2}).
\end{equation}         
and
\begin{equation}\label{w}
w =  - z{u_{0,x}} + \mu \frac{{{z^3}}}{6}{u_{0,xxx}} + O({\mu ^2}).    
\end{equation}
Substituting the expressions of $u,w$ in (\ref{finalu})-(\ref{w}) into the second equation in Eq. (\ref{scaling and shear}), we find
\begin{equation*}
{p_z} = \mu z({u_{0,tx}} + \varepsilon {u_0}{u_{0,xx}} -\varepsilon u_{0,x}^2) + \mu A{z^2}{u_{0,xx}}+O({\mu ^2}).
\end{equation*}
Thanks to the boundary condition of $p$ on $z=h$, integrating the above equality with respect to $z$, we get
\[p = \eta  - \frac{\mu }{2}({h^2} - {z^2})({u_{0,tx}} + \varepsilon {u_0}{u_{0,xx}} - \varepsilon u_{0,x}^2) - \mu A\frac{{{h^3} - {z^3}}}{3}{u_{0,xx}} + O({\mu ^2}).\]
Differentiating the above equation with respect to $x$, we obtain
\begin{equation*}
\begin{split}
{p_x} =&{\eta _x} - \frac{\mu }{2}({h^2} - {z^2}){({u_{0,tx}} + \varepsilon {u_0}{u_{0,xx}} - \varepsilon u_{0,x}^2)_x} - \mu h{h_x}({u_{0,tx}}
+ \varepsilon {u_0}{u_{0,xx}} - \varepsilon u_{0,x}^2)\\
&- \mu A\frac{{{h^3} - {z^3}}}{3}{u_{0,xxx}} - \mu A{h^2}{h_x}{u_{0,xx}} + O({\mu ^2}).
\end{split}
\end{equation*}
Plugging the expressions of $u,w$ in (\ref{finalu})-(\ref{w}) into the first equation in Eq. (\ref{scaling and shear}), it yields
\begin{equation*}
\begin{split}
 -{p_x} = {u_{0,t}}& - \mu \frac{{{z^2}}}{2}{u_{0,txx}} + Az({u_{0,x}} - \mu \frac{{{z^2}}}{2}{u_{0,xxx}})- Az{u_{0,x}} + \mu A\frac{{{z^3}}}{6}{u_{0,xxx}} + \varepsilon {u_0}{u_{0,x}}\\
  &+ \varepsilon \mu \frac{{{z^2}}}{2}({u_{0,x}}{u_{0,xx}} - {u_0}{u_{0,xxx}}) + O({\mu ^2}).
\end{split}
\end{equation*}
It thus follows the last two equations that
\begin{equation}\label{u0}
\begin{split}
{u_{0,t}}&+ \varepsilon {u_0}{u_{0,x}} + {\eta _x}= \frac{\mu }{2}{h^2}{({u_{0,tx}} + \varepsilon {u_0}{u_{0,xx}} - \varepsilon u_{0,x}^2)_x}+ \mu A\frac{{{h^3}}}{3}{u_{0,xxx}}\\
&+ \mu h{h_x}({u_{0,tx}} + \varepsilon {u_0}{u_{0,xx}} - \varepsilon u_{0,x}^2)
 + \mu A{h^2}{h_x}{u_{0,xx}} + O({\mu ^2}).
\end{split}
\end{equation}
Secondly, we apply the following relation between $\bar{u}$ and $u_0$ to find the second equation of the Green-Naghdi equations with a linear shear.
\begin{equation}\label{relation}
 \bar{u } = {u_0} - \mu\frac{{{h^2}}}{6}{u_{0,xx}} + O({\mu ^2}),
\end{equation}
which can be obtained by combining (\ref{ubar}) with (\ref{finalu}). Then, we deduce the following more facts about $\bar{u}$ and $u_0$ from (\ref{relation}).
\begin{equation}\label{relation1}
\mu  \bar{u}  = \mu {u_0} + O({\mu ^2}),
\end{equation}
\begin{equation}\label{relation2}
{ \bar{u} _t} = {u_{0,t}} - \mu \frac{{{h^2}}}{6}{u_{0,txx}} - \mu \frac{{h{h_t}}}{3}{u_{0,xx}} + O({\mu ^2}),
\end{equation}
and
\begin{equation}\label{relation3}
\frac{\mu }{{3h}}({({h^3}({{\bar u}_{xt}} + \varepsilon \bar u{{\bar u}_{xx}} - \varepsilon \bar u_x^2))_x}
= \frac{\mu }{{3h}}{({h^3}({u_{0,xt}} + \varepsilon {u_0}{u_{0,xx}} - \varepsilon u_{0,x}^2))_x} + O({\mu ^2}).
\end{equation}
Now we claim that
\begin{equation}\label{eq}
\begin{split}
\varepsilon \bar u{{\bar u}_x} &= \varepsilon {u_0}{u_{0,x}} - \varepsilon \mu \frac{{{h^2}}}{6}{u_0}{u_{0,xxx}} + \mu \frac{{h{h_t}}}{3}{u_{0,xx}} + \mu \frac{A}{3}{h^2}{h_x}{u_{0,xx}} + \varepsilon \mu \frac{{{h^2}}}{6}{u_{0,x}}{u_{0,xx}} + O({\mu ^2}).
\end{split}
\end{equation}
To prove this, in view of the expression of $\bar u$ in (\ref{relation}), we have
\begin{equation}\label{eq1}
\varepsilon \bar u{{\bar u}_x} = \varepsilon {u_0}{u_{0,x}} - \varepsilon \mu \frac{{{h^2}}}{6}{u_0}{u_{0,xxx}} - \varepsilon \mu \frac{{h{h_x}}}{3}{u_0}{u_{0,xx}} - \varepsilon \mu \frac{{{h^2}}}{6}{u_{0,x}}{u_{0,xx}} + O({\mu ^2}).\end{equation}
According to (\ref{h}), we find
\begin{equation}\label{eq2}
\begin{split}
\mu \frac{{h{h_t}}}{3}{u_{0,xx}}
&= \frac{{\mu h}}{3}{u_{0,xx}} \cdot {h_t}
                = \frac{{\mu h}}{3}{u_{0,xx}}( - Ah{h_x} - \varepsilon {h_x}\bar u - \varepsilon h{{\bar u}_x}) \\
                & =  - \frac{\mu }{3}A{h^2}{h_x}{u_{0,xx}} - \frac{{\varepsilon \mu }}{3}h{h_x}{u_0}{u_{0,xx}} - \frac{{\varepsilon \mu }}{3}{h^2}{u_{0,x}}{u_{0,xx}}+ O({\mu ^2}).
\end{split}
\end{equation}
Hence we obtain (\ref{eq}) by plugging (\ref{eq2}) into the expression of $\varepsilon \bar u{{\bar u}_x}$ in (\ref{eq1}). Now subtracting $\frac{{\mu {h^2}}}{6}{({u_{0,tx}} + \varepsilon {u_0}{u_{0,xx}} - \varepsilon u_{0,x}^2)_x}$ on both sides of (\ref{u0}), we get
\begin{equation*}
\begin{split}
{u_{0,t}}-& \frac{\mu }{6}{h^2}{({u_{0,tx}} + \varepsilon {u_0}{u_{0,xx}} - \varepsilon u_{0,x}^2)_x} + \varepsilon {u_0}{u_{0,x}} + {\eta _x}\\
=& \frac{\mu }{2}{h^2}{({u_{0,tx}} + \varepsilon {u_0}{u_{0,xx}} - \varepsilon u_{0,x}^2)_x} - \frac{\mu }{6}{h^2}{({u_{0,tx}}
+ \varepsilon {u_0}{u_{0,xx}} - \varepsilon u_{0,x}^2)_x}\\
 &+ \mu h{h_x}({u_{0,tx}} + \varepsilon {u_0}{u_{0,xx}} - \varepsilon u_{0,x}^2) + \mu A\frac{{{h^3}}}{3}{u_{0,xxx}} + \mu A{h^2}{h_x}{u_{0,xx}} + O({\mu ^2}),
\end{split}
\end{equation*}
which is equivalent to
\begin{equation}\label{finalu0}
\begin{split}
{u_{0,t}}&- \frac{\mu }{6}{h^2}{u_{0,txx}} - \mu \frac{{h{h_t}}}{3}{u_{0,xx}}+\varepsilon {u_0}{u_{0,x}} + \mu \frac{{h{h_t}}}{3}{u_{0,xx}} - \frac{{\mu {h^2}}}{6}{(\varepsilon {u_0}{u_{0,xx}} - \varepsilon u_{0,x}^2)_x}+ {\eta _x}\\
      & +\mu \frac{A}{3}{h^2}{h_x}{u_{0,xx}}= \frac{\mu }{3}{h^2}{({u_{0,tx}} + \varepsilon {u_0}{u_{0,xx}} - \varepsilon u_{0,x}^2)_x} + \mu h{h_x}({u_{0,tx}} + \varepsilon {u_0}{u_{0,xx}} - \varepsilon u_{0,x}^2)\\
 &+ \mu A\frac{{{h^3}}}{3}{u_{0,xxx}} + \mu A{h^2}{h_x}{u_{0,xx}}+\mu\frac{A}{3}{h^2}{h_x}{u_{0,xx}} + O({\mu ^2}).
\end{split}
\end{equation}
In light of (\ref{relation1})-(\ref{eq}), we finally establish the second equation of the Green-Naghdi equations with a linear shear from (\ref{finalu0}).
\begin{equation}\label{finalubar}
\begin{split}
{{\bar u}_t} + \varepsilon \bar u{{\bar u}_x} + {\eta _x}
=& \frac{\mu }{{3(1 + \varepsilon \eta )}}{({(1 + \varepsilon \eta )^3}({{\bar u}_{xt}} + \varepsilon \bar u{{\bar u}_{xx}} - \varepsilon \bar u_x^2))_x}\\&+ \mu A\frac{{{h^3}}}{3}{{\bar u}_{xxx}} + \mu A\frac{4}{3}{h^2}{h_x}{{\bar u}_{xx}}+O({\mu ^2}).
\end{split}
\end{equation}
For simplicity, dropping the superscript bar in $\bar u$ of (\ref{eta}) and (\ref{finalubar}) in the sequel, we obtain
\begin {equation}
\left\{
\begin{aligned}\label{GNwithshear}
&{{\eta _t} + A(1 + \varepsilon \eta ){\eta _x} + {{((1 + \varepsilon \eta )u)}_x} = 0}, \\
&{u_t} + \varepsilon u{u_x} + {\eta _x}= \frac{\mu }{{3(1 + \varepsilon \eta )}}{{({{(1 + \varepsilon \eta )}^3}({u_{xt}} + \varepsilon u{u_{xx}} - \varepsilon u_x^2))}_x}\\
&\quad \quad \quad \quad \quad \quad \quad \quad +\mu A\frac{{{{(1 + \varepsilon \eta )}^3}}}{3}{u_{xxx}} + \varepsilon \mu A\frac{{4{{(1 + \varepsilon \eta )}^2}}}{3}{\eta _x}{u_{xx}}+O({\mu ^2}).
\end{aligned}
\right.
\end{equation}

\subsection{Two-component system with constant vorticity under CH scaling}
In this subsection, inspired by the work of Constantin and Ivanov \cite{Constantin-I},
we start from Eq. (\ref{GNwithshear}) to proceed with the derivation of our two-component system under CH scaling ($\mu\ll1, \varepsilon=O(\sqrt{\mu})$).

Up to the correction $O(\varepsilon,\mu),$ we deduce from Eq. (\ref{GNwithshear}) the linearized system.
\begin{equation}\label{linearueta}
\left\{
\begin{aligned}
&{\eta _t}+ A{\eta _x}+u_x =O(\varepsilon,\mu),\\
&{u_t}+{\eta _x}=O(\varepsilon,\mu),
\end{aligned} \right.
\end{equation}
which implies
\begin{equation*}
\left\{
\begin{aligned}
{\eta _{tt}} + A{\eta _{tx}} - {\eta _{xx}} = O(\varepsilon ,\mu ),\\
{u_{tt}} + A{u_{tx}} - {u_{xx}} = O(\varepsilon ,\mu ).
\end{aligned}
\right.
\end{equation*}
Then we get the general solution of the above system.
\begin{equation*}
\left\{
\begin{aligned}
\eta  = \eta_1 (x - \frac{\sqrt{A^{2}+4}+A}{2}t) + \eta_2 (x + \frac{\sqrt{A^{2}+4}-A}{2}t) + O(\varepsilon ,\mu ),\\
u = u_1(x - \frac{\sqrt{A^{2}+4}+A}{2}t) + u_2 (x + \frac{\sqrt{A^{2}+4}-A}{2}t) + O(\varepsilon ,\mu ).
\end{aligned}
\right.
\end{equation*}
It is convenient first to restrict ourselves to waves moving towards the right side, i.e.,
\begin{equation*}
\left\{
\begin{aligned}
\eta  = \eta (x - ct) + O(\varepsilon ,\mu ),\\
u = u(x - ct) + O(\varepsilon ,\mu ),
\end{aligned}
\right.
\end{equation*}
where $c = \frac{\sqrt{A^{2}+4}+A}{2}$, which have been confirmed by the Burns condition. Thus we have
\begin{equation}\label{linearueta1}
\left\{
 \begin{aligned}
 {\eta _t} =  - c{\eta _x} + O(\varepsilon ,\mu ),\\
{u_t} =  - c{u_x} + O(\varepsilon ,\mu ).
\end{aligned}
\right.
\end{equation}
It then follows from (\ref{linearueta}) and (\ref{linearueta1}) that
\begin{equation}\label{relationueta}
{u_t} + {\eta _x} = O(\varepsilon ,\mu ),\quad\quad       u - \frac{1}{c}\eta  = O(\varepsilon ,\mu ),\quad\quad
{u_x} + \frac{1}{{{c^2}}}{\eta _t} = O(\varepsilon ,\mu ).
\end{equation}

Under the CH scaling, we expand Eq. (\ref{GNwithshear}) up to order of $O(\varepsilon^2 \mu)$. It reads,
\begin{equation}\label{GN2}
\left\{
\begin{aligned}
&{\eta _t} + A(1 + \varepsilon \eta ){\eta _x} + {((1 + \varepsilon \eta )u)_x} = 0,\\
&{(u - \frac{\mu }{3}{u_{xx}})_t} + {\eta _x} + \varepsilon u{u_x} =  - \varepsilon \mu {c^2}{u_x}{u_{xx}} - \frac{{\varepsilon \mu }}{3}{u_x}{u_{xx}} + \frac{{\varepsilon \mu }}{3}u{u_{xxx}}\\
&\qquad\qquad- \frac{{2\varepsilon \mu }}{3}{c^2}u{u_{xxx}}+ \varepsilon \mu Acu{u_{xxx}} + \frac{{4\varepsilon \mu }}{3}Ac{u_x}{u_{xx}} + \frac{\mu }{3}A{u_{xxx}} + O({\varepsilon ^2}\mu ,{\mu ^2}).
\end{aligned}
\right.
\end{equation}
To obtain our system, we introduce an auxiliary variable as in \cite{Constantin-I}.
\begin{equation}\label{rho}
\rho  = 1 + \frac{1}{2}\varepsilon \eta  - \frac{1}{8}{\varepsilon ^2}({u^2} + {\eta ^2}),
\end{equation}
which gives
\begin{equation}\label{rhorho}
{\rho ^2} = 1 + \varepsilon \eta  - \frac{1}{4}{\varepsilon ^2}{u^2} - \frac{1}{8}{\varepsilon ^3}({u^2} + {\eta ^2})\eta  + O({\varepsilon ^4}).
\end{equation}
In view of (\ref{rhorho}), $\eta_x$ can be represented by
\begin{equation}\label{etax}
{\eta _x} = \frac{1}{\varepsilon }({\rho ^2} + \frac{1}{4}{\varepsilon ^2}{u^2} + \frac{1}{8}{\varepsilon ^3}({u^2} + {\eta ^2})\eta )_x + O({\varepsilon ^3}).
\end{equation}
Combining (\ref{relationueta})-(\ref{GN2}) with (\ref{rho}), we obtain
\begin{equation}\label{rhot}
\begin{aligned}
{\rho _t}&= \frac{1}{2}\varepsilon {\eta _t} - \frac{1}{8}{\varepsilon ^2}{({u^2} + {\eta ^2})_t}\vspace{1ex}= - \frac{1}{2}\varepsilon [A(1 + \varepsilon \eta ){\eta _x} + {((1 + \varepsilon \eta )u)_x}] - \frac{1}{4}{\varepsilon ^2}(u{u_t} + \eta {\eta _t})\vspace{1ex}\\
 &=  - \frac{1}{2}\varepsilon A(1 + \varepsilon \eta ){\eta _x} - \frac{1}{2}\varepsilon {((1 + \varepsilon \eta )u)_x} + \frac{1}{4}{\varepsilon ^2}u({\eta _x} + \varepsilon u{u_x})\\
  &\quad + \frac{1}{4}{\varepsilon ^2}\eta [A(1 + \varepsilon \eta ){\eta _x} + {((1 + \varepsilon \eta )u)_x}]\vspace{1ex}+O({\varepsilon ^2}\mu ,{\varepsilon ^4})\\
 &=  - \frac{1}{2}\varepsilon A(1 + \varepsilon \eta ){\eta _x} - \frac{1}{2}\varepsilon {((1 + \varepsilon \eta )u)_x} + \frac{1}{4}{\varepsilon ^2}u{\eta _x} + \frac{1}{4}{\varepsilon ^3}{u^2}{u_x} + \frac{1}{4}{\varepsilon ^2}A(1 + \varepsilon \eta )\eta {\eta _x} \\
 &\quad+ \frac{1}{4}{\varepsilon ^2}\eta {u_x} + \frac{1}{4}{\varepsilon ^3}\eta {(\eta u)_x}\vspace{1ex}+O({\varepsilon ^2}\mu ,{\varepsilon ^4})\\
 &=  - \frac{1}{2}\varepsilon A(1 + \varepsilon \eta ){\eta _x} - \frac{1}{2}\varepsilon {((1 + \frac{1}{2}\varepsilon \eta )u)_x} + \frac{1}{4}{\varepsilon ^3}{u^2}{u_x}+ \frac{1}{4}{\varepsilon ^2}A(1 + \varepsilon \eta )\eta {\eta _x} \\
 &\quad + \frac{1}{4}{\varepsilon ^3}\eta {(\eta u)_x}\vspace{1ex}+O({\varepsilon ^2}\mu ,{\varepsilon ^4})\\
 & =  - \frac{1}{2}\varepsilon {(\rho u)_x} - \frac{1}{{16}}{\varepsilon ^3}{({u^3} + {\eta ^2}u)_x} + \frac{1}{{12}}{\varepsilon ^3}{({u^3})_x} + \frac{{{c^2}}}{6}{\varepsilon ^3}{({u^3})_x}\vspace{1ex} - \frac{1}{2}\varepsilon A(1 + \varepsilon \eta ){\eta _x} \\
 &\quad+ \frac{1}{4}{\varepsilon ^2}A(1 + \varepsilon \eta )\eta {\eta _x}\vspace{1ex}+O({\varepsilon ^2}\mu ,{\varepsilon ^4})\\
 &=  - \frac{1}{2}\varepsilon {(\rho u)_x} - \frac{{{c^2} + 1}}{{16}}{\varepsilon ^3}{({u^3})_x} + \frac{1}{{12}}{\varepsilon ^3}{({u^3})_x} + \frac{{{c^2}}}{6}{\varepsilon ^3}{({u^3})_x}\vspace{1ex} - \frac{1}{2}\varepsilon A(1 + \varepsilon \eta ){\eta _x}\\
 &\quad + \frac{1}{4}{\varepsilon ^2}A(1 + \varepsilon \eta )\eta {\eta _x}+O({\varepsilon ^2}\mu ,{\varepsilon ^4}).
\end{aligned}
\end{equation}
For the latter two terms in the right hand side of (\ref{rhot}), we deduce from (\ref{etax}) that
\begin{equation}\label{twoterms}
\begin{aligned}
 &- \frac{1}{2}\varepsilon A(1 + \varepsilon \eta ){\eta _x} + \frac{1}{4}{\varepsilon ^2}A(1 + \varepsilon \eta )\eta {\eta _x}\vspace{1ex}=  - \frac{1}{2}\varepsilon A{\eta _x} - \frac{1}{2}{\varepsilon ^2}A\eta {\eta _x} + \frac{1}{4}{\varepsilon ^2}A\eta {\eta _x} + \frac{1}{4}{\varepsilon ^3}A{\eta ^2}{\eta _x}\vspace{1ex}\\
 &=  - \frac{1}{2}\varepsilon A{\eta _x} - \frac{1}{4}{\varepsilon ^2}A\eta {\eta _x} + \frac{1}{4}{\varepsilon ^3}A{\eta ^2}{\eta _x}\vspace{1ex}=  - \frac{1}{2}\varepsilon A{\eta _x}(1 + \frac{1}{2}\varepsilon \eta ) + \frac{1}{4}{\varepsilon ^3}A{\eta ^2}{\eta _x}\vspace{1ex}\\
 &=  - \frac{1}{2}\varepsilon A{\eta _x}[\rho  + \frac{1}{8}{\varepsilon ^2}({u^2} + {\eta ^2})] + \frac{1}{4}{\varepsilon ^3}A{\eta ^2}{\eta _x}\vspace{1ex}=  - \frac{1}{2}\varepsilon A\rho {\eta _x} - \frac{1}{{16}}{\varepsilon ^3}A({u^2} + {\eta ^2}){\eta _x} + \frac{1}{4}{\varepsilon ^3}A{\eta ^2}{\eta _x}\vspace{1ex}\\
 &=  - \frac{1}{2}A\rho {({\rho ^2} + \frac{1}{4}{\varepsilon ^2}{u^2} + \frac{1}{8}{\varepsilon ^3}({u^2} + {\eta ^2})\eta )_x} - \frac{{{c^3} + c}}{{48}}{\varepsilon ^3}A{({u^3})_x} + \frac{{{c^3}}}{{12}}{\varepsilon ^3}A{({u^3})_x} + O({\varepsilon ^2}\mu ,{\mu ^2},{\varepsilon ^4})\vspace{1ex}\\
 &=  - \frac{1}{3}A{({\rho ^3})_x} - \frac{1}{8}{\varepsilon ^2}A\rho {({u^2})_x} - \frac{{{c^3} + c}}{{16}}{\varepsilon ^3}A\rho {({u^3})_x} + \frac{{3{c^3} - c}}{{48}}{\varepsilon ^3}A{({u^3})_x} + O({\varepsilon ^2}\mu ,{\mu ^2},{\varepsilon ^4}).
\end{aligned}
\end{equation}
Plugging (\ref{twoterms}) into (\ref{rhot}), we find
\begin{equation*}
\begin{aligned}
{\rho _t} &= - \frac{1}{3}A{({\rho ^3})_x} - \frac{1}{2}\varepsilon {(\rho u)_x} - \frac{1}{8}{\varepsilon ^2}A\rho {({u^2})_x} + (\frac{{5{c^2} + 1}}{{48}} + \frac{({3{c^3} - c})A}{{48}}){\varepsilon ^3}{({u^3})_x} \vspace{1ex}\\
&\quad\quad - \frac{{{c^3} + c}}{{16}}{\varepsilon ^3}A\rho {({u^3})_x}+  O({\varepsilon ^2}\mu ,{\mu ^2},{\varepsilon ^4}).
\end{aligned}
\end{equation*}
Noticing that ${\varepsilon ^3} {({u^3})_x} = {\varepsilon ^3}\rho{({u^3})_x} + O({\varepsilon ^4})$, we have
\begin{equation}\label{equationrho}
{\rho _t} + \frac{1}{3}A{({\rho ^3})_x} + \frac{1}{2}\varepsilon {(\rho u)_x} + \frac{1}{8}{\varepsilon ^2}A\rho {({u^2})_x} + {\beta _1}{\varepsilon ^3}\rho {({u^3})_x} =  O({\varepsilon ^2}\mu ,{\mu ^2},{\varepsilon ^4}),
\end{equation}
where ${\beta _1} =  - (\frac{{5{c^2} + 1}}{{48}} + \frac{({3{c^3} - c})A}{{48}} - \frac{{A({c^3} + c)}}{{16}}) =  - \frac{{{c^2} + 5}}{{48}}$.

On the other hand, denoting $m=u - \frac{\mu }{3}{u_{xx}}$, it follows from (\ref{etax}) that
\[{m_t} + \varepsilon u{u_x} + \frac{1}{\varepsilon }({\rho ^2} + \frac{1}{4}{\varepsilon ^2}{u^2} + \frac{1}{8}{\varepsilon ^3}({u^2} + {\eta ^2})\eta)_x
= \frac{\mu }{3}A{u_{xxx}} + O(\varepsilon \mu ,{\varepsilon ^3}).\]
On account of (\ref{relationueta}), we get
\[\varepsilon {m_t} + \frac{3}{2}{\varepsilon ^2}u{u_x} + {({\rho ^2})_x} + {\alpha _1}{\varepsilon ^3}{({u^3})_x} - \frac{{\varepsilon \mu }}{3}A{u_{xxx}} = O({\varepsilon ^2}\mu ,{\mu ^2},{\varepsilon ^4}).\]
where ${\alpha _1} = \frac{{{c^3} + c}}{8}$. Moreover, using (\ref{relationueta}) and (\ref{rho}), we have
\begin{equation*}
\begin{aligned}
{\alpha _1}{\varepsilon ^3}{({u^3})_x} &= {\alpha _2}{\varepsilon ^3}{({u^3})_x} + ({\alpha _1} - {\alpha _2}){\varepsilon ^3}{({u^3})_x}\vspace{1ex}= {\alpha _2}{\varepsilon ^3}{({u^3}\rho )_x} + ({\alpha _1} - {\alpha _2}){\varepsilon ^3}{({u^3})_x} + O({\varepsilon ^4})\vspace{1ex}\\
                &= \frac{{2{\alpha _2}}}{c}{\varepsilon ^2}{({u^2}\rho  \cdot \frac{{\varepsilon \eta }}{2})_x} + ({\alpha _1} - {\alpha _2}){\varepsilon ^3}{({u^3})_x} + O({\varepsilon ^4})\vspace{1ex}\\
                &= \frac{{2{\alpha _2}}}{c}{\varepsilon ^2}{({u^2}\rho (\rho  - 1))_x} + ({\alpha _1} - {\alpha _2}){\varepsilon ^3}{({u^3})_x} + O({\varepsilon ^4}).
\end{aligned}
\end{equation*}
where $\alpha _2$ is a constant yet to be determined. With the choice $6{\beta _1}=\frac{\alpha _2}{c}$ (thus, ${\alpha _2} =  - \frac{{c({c^2} + 5)}}{8}$), we obtain
\begin{equation}\label{equationu}
\begin{split}
\varepsilon {m_t} + \frac{3}{2}{\varepsilon ^2}u{u_x} + {({\rho ^2})_x} + 12{\beta _1}{\varepsilon ^2}{({u^2}\rho (\rho  - 1))_x}+ {\beta _2}{\varepsilon ^3}{({u^3})_x} - \frac{{\varepsilon \mu }}{3}A{u_{xxx}} = O({\varepsilon ^2}\mu ,{\mu ^2},{\varepsilon ^4}),
 \end{split}
\end{equation}
where ${\beta _2} = {\alpha _1} - {\alpha _2} = \frac{{{c^3} + 3c}}{4}$. For some constant $\sigma$, using the fact
\[\frac{3}{2}{\varepsilon ^2}u{u_x} = \frac{{\sigma {\varepsilon ^2}}}{2}(2m{u_x} + u{m_x}) + \frac{3}{2}(1 - \sigma ){\varepsilon ^2}u{u_x} + O({\varepsilon ^2}\mu ),\]
we can rewrite Eq.(\ref{equationu}) as
\begin{equation}\label{equationu1}
\begin{aligned}
{\varepsilon {m_t} }&+ \frac{{\sigma {\varepsilon ^2}}}{2}(2m{u_x} + u{m_x}) + \frac{3}{2}(1 - \sigma ){\varepsilon ^2}u{u_x} + {({\rho ^2})_x}\\
&\quad + 12{\beta _1}{\varepsilon ^2}{({u^2}\rho (\rho  - 1))_x} + {\beta _2}{\varepsilon ^3}{({u^3})_x} - \frac{{\varepsilon \mu }}{3}A{u_{xxx}} = O({\varepsilon ^2}\mu ,{\mu ^2},{\varepsilon ^4}).
\end{aligned}
\end{equation}
Then applying the rescaling $u \to \frac{2}{\varepsilon }u,\ x \to \sqrt {\frac{\mu }{3}}x ,\ t \to \sqrt {\frac{\mu }{3}}t $ to (\ref{equationrho}) and (\ref{equationu1}), we
obtain the following two-component system with constant vorticity under the CH scaling
\begin{equation*}
\left\{
 \begin{aligned}
&{\rho _t} + \frac{1}{3}A{({\rho ^3})_x} + {(\rho u)_x} + \frac{1}{2}A\rho {({u^2})_x} + 8{\beta _1}\rho {({u^3})_x} =  0,\vspace{1ex}\\
&{m_t} + \sigma (2m{u_x} + u{m_x}) + 3(1 - \sigma )u{u_x} + \frac{1}{2}{({\rho ^2})_x}\vspace{1ex}\\
 &\qquad\qquad\qquad + 24{\beta _1}{({u^2}\rho (\rho  - 1))_x} + 4{\beta _2}{({u^3})_x} - A{u_{xxx}} = 0,
\end{aligned}
\right.
\end{equation*}
where $m = u - {u_{xx}},$ and the coefficients ${\beta _1}= - \frac{{{c^2} + 5}}{{48}},$ ${\beta _2}= \frac{{{c^3} + 3c}}{4},$ $A= \frac{{{c^2} - 1}}{c}.$

\section{Local well-posedness}\label{sec3}
\newtheorem{theorem3}{Theorem}[section]
\newtheorem{lemma3}{Lemma}[section]
\newtheorem {remark3}{Remark}[section]
\newtheorem {definition3}{Definition}[section]
\newtheorem{corollary3}{Corollary}[section]
\par

In this section, we shall establish the local well-posedness of the Cauchy problem of the general system (\ref{two-component}) in nonhomogeneous Besov spaces.
Firstly, we obtain the uniqueness and continuity with respect to the initial data in some appropriate norms by the following a priori estimates.
\begin{lemma3}\label{Lem3.1}
Let $1\leq p,r\leq \infty$ and $s>2 + \frac{1}{p}$. Suppose that $(u^{(i)},\zeta^{(i)})\in {L^\infty }(0,T;B_{p,r}^s\times B_{p,r}^{s-1}) \cap C([0,T];\mathcal{S'} \times \mathcal{S'})$, $(i=1,2)$
be two given solutions of the system $(\ref{two-component})$ with the initial data $(u_0^{(i)},\zeta_0^{(i)})\in B_{p,r}^s\times B_{p,r}^{s-1}$, $(i=1,2)$. Denoting
$u^{(12)}:=u^{(2)}-u^{(1)}$ and $\zeta^{(12)} :=\zeta^{(2)}-\zeta^{(1)}$, then for all $t\in[0,T]$, we have
\begin{eqnarray}\label{uniqueness}
{\| u^{(12)} \|_{B_{p,r}^{s - 1}}}+ {\| \zeta^{(12)} \|_{B_{p,r}^{s - 2}}}&\le& ({\| u_0^{(12)}\|_{B_{p,r}^{s - 1}}} + {\| \zeta_0^{(12)}\|_{B_{p,r}^{s - 2}}})\nonumber \\
 &&\times e^{ C\int_0^t \big( 1+\sum\limits_{k = 1}^3 ({{{\| {{u^{(1)}}} \|}_{B_{p,r}^s}} + {{\| {{u^{(2)}}} \|}_{B_{p,r}^s}} + {{\| {{\zeta ^{(1)}}} \|}_{B_{p,r}^{s - 1}}} + {{\| {{\zeta ^{(2)}}} \|}_{B_{p,r}^{s - 1}}}){^k}}\big) d\tau}.\nonumber \\
\end{eqnarray}
\end{lemma3}
\begin{proof}
It is obvious that $(u^{(12)},\zeta^{(12)})\in {L^\infty }(0,T;B_{p,r}^s\times B_{p,r}^{s-1}) \cap C([0,T];\mathcal{S'} \times \mathcal{S'})$ and solves the following Cauchy problem of the transport
equations
\begin{equation}\label{transport}
\left\{
{\begin{aligned}
&{{\partial _t}{u^{(12)}} + \sigma {u^{(2)}}{\partial _x}{u^{(12)}} = {F}(t,x)},\\
&{{\partial _t}{\zeta ^{(12)}} + ({u^{(2)}} + {b_1}{{({\zeta ^{(2)}})}^2} + 2{b_1}{\zeta ^{(2)}}){\partial _x}{\zeta ^{(12)}} =  {G}(t,x)},\\
&{{u^{(12)}}(0,x)= u_0^{12}},\\
&{{\zeta ^{(12)}}(0,x) = \zeta _0^{12}},
\end{aligned}}
\right.
\end{equation}
where
\begin{equation*}
\begin{aligned}
{F} (t,x)&:=  - \sigma {u^{(12)}}{\partial _x}{u^{(1)}} +P(D)\big(\frac{{3 - \sigma }}{2}{u^{(12)}}({u^{(1)}} + {u^{(2)}}) + \frac{\sigma }{2}{\partial _x}{u^{(12)}}{\partial _x}({u^{(1)}} + {u^{(2)}})\\
 &+ \frac{1}{2}{\zeta ^{(12)}}({\zeta ^{(1)}} + {\zeta ^{(2)}} + 2)+ {a_1}( {({\zeta ^{(2)}})^2}{u^{(12)}} + {\zeta ^{(12)}}({\zeta ^{(1)}} + {\zeta ^{(2)}}){u^{(1)}})  \\
 &+ (2{a_1} + {a_2})( {\zeta ^{(2)}}{u^{(12)}} + {\zeta ^{(12)}}{u^{(1)}})+ ({a_1} + {a_2} + {a_3} + {a_4}){u^{(12)}}\\
 & + {a_5}( {({u^{(2)}})^2}{\zeta ^{(12)}}({\zeta ^{(1)}} + {\zeta ^{(2)}} + 1))  + {a_5}( {u^{(12)}}({u^{(1)}} + {u^{(2)}}){({\zeta ^{(1)}})^2})\\
 & + {a_5}( {u^{(12)}}({u^{(1)}} + {u^{(2)}}){\zeta ^{(1)}})+ {a_6}( {u^{(12)}}({({u^{(2)}})^2} + {u^{(1)}}{u^{(2)}} + {({u^{(1)}})^2})) \big),
\end{aligned}
\end{equation*}
and
\begin{equation*}
\begin{aligned}
{G}(t,x) &:=  - {u^{(12)}}{\partial _x}{\zeta ^{(1)}} - {b_1}({\zeta ^{(1)}} + {\zeta ^{(2)}}){\zeta ^{(12)}}{\partial _x}{\zeta ^{(1)}} - 2{b_1}{\zeta ^{(12)}}{\partial _x}{\zeta ^{(1)}} - {\zeta ^{(12)}}{\partial _x}{u^{(2)}}\\
 &- {\zeta ^1}{\partial _x}{u^{(12)}} - {\partial _x}{u^{(12)}} - \frac{{{b_2}}}{2}{\zeta ^{(12)}}{\partial _x}({u^{(2)}})^2- {b_2}{\zeta ^{(1)}}({u^{(2)}}{\partial _x}{u^{(12)}} + {u^{(12)}}{\partial _x}{u^{(1)}})\\
  &- {b_2}({u^{(2)}}{\partial _x}{u^{(12)}} + {u^{(12)}}{\partial _x}{u^{(1)}})- \frac{{{b_3}}}{3}{\zeta ^{(12)}}{\partial _x}({u^{(2)}})^3 - {b_3}{({u^{(2)}})^2}{\partial _x}{u^{(12)}} \\
 &- {b_3}{\zeta ^{(1)}}({({u^{(2)}})^2}{\partial _x}{u^{(12)}} +{u^{(12)}}({u^{(1)}} + {u^{(2)}}){\partial _x}{u^{(1)}}) -b_3 {u^{(12)}}({u^{(1)}} +{u^{(2)}}){\partial _x}{u^{(1)}}.
\end{aligned}
\end{equation*}
According to Lemma \ref{Lem6.1}, $B_{p,r}^{s-1}, B_{p,r}^{s-2} $ with $s>2 +\frac{1}{p}$ are both algebras, and $P(D)\in Op(S^{-1})$, we have
\begin{eqnarray*}
{\| {{u^{(12)}}{\partial _x}{u^{(1)}}} \|_{B_{p,r}^{s - 1}}} \le C{\| {{u^{(12)}}} \|_{B_{p,r}^{s - 1}}}{\| {{u^{(1)}}} \|_{B_{p,r}^s}},
\end{eqnarray*}
\begin{eqnarray*}
&{\| {P(D)\big(
{u^{(12)}}(\frac{{3 - \sigma }}{2}({u^{(1)}} + {u^{(2)}})+({a_1} + {a_2} + {a_3} + {a_4})+a_6({({u^{(2)}})^2} + {u^{(1)}}{u^{(2)}} + {({u^{(1)}})^2})
)\big)\|_{B_{p,r}^{s - 1}}}}\\
&\le C{\| {{u^{(12)}}} \|_{B_{p,r}^{s - 1}}}\big(1+({\| {{u^{(1)}}} \|_{B_{p,r}^s}} + {\| {{u^{(2)}}} \|_{B_{p,r}^s}})+
({\| {{u^{(1)}}} \|_{B_{p,r}^s}} + {\| {{u^{(2)}}} \|_{B_{p,r}^s}})^2\big),
\end{eqnarray*}
\begin{eqnarray*}
&{\| {P(D)\big({u^{(12)}}({a_1}({\zeta ^{(2)}})^2+(2{a_1} + {a_2}) {\zeta ^{(2)}}+{a_5}({u^{(1)}} + {u^{(2)}}){({\zeta ^{(1)}})^2}+
{a_5}({u^{(1)}} + {u^{(2)}}){\zeta ^{(1)}})\big)\|_{B_{p,r}^{s - 1}}}}\\
&\le C{\| {{u^{(12)}}} \|_{B_{p,r}^{s - 1}}}\big(({\| {{\zeta^{(2)}}} \|_{B_{p,r}^s}} + {\| {{\zeta^{(2)}}} \|^2_{B_{p,r}^s}})
+({\| {{u^{(1)}}} \|_{B_{p,r}^s}} + {\| {{u^{(2)}}} \|_{B_{p,r}^s}})({\| {{\zeta^{(1)}}} \|_{B_{p,r}^s}} + {\| {{\zeta^{(1)}}} \|^2_{B_{p,r}^s}})
\big),
\end{eqnarray*}
\begin{eqnarray*}
{\| {P(D)({\partial _x}{u^{(12)}}{\partial _x}({u^{(1)}} + {u^{(2)}})
)\|_{B_{p,r}^{s - 1}}}}\le C{\| {{u^{(12)}}} \|_{B_{p,r}^{s - 1}}}({\| {{u^{(1)}}} \|_{B_{p,r}^s}} + {\| {{u^{(2)}}} \|_{B_{p,r}^s}}),
\end{eqnarray*}
\begin{eqnarray*}
{\| {P(D)({\zeta ^{(12)}}({\zeta ^{(1)}} + {\zeta ^{(2)}} + 2))\|_{B_{p,r}^{s - 1}}}}
\le C {\| {{\zeta^{(12)}}} \|_{B_{p,r}^{s - 2}}}(1+{\| {{\zeta^{(1)}}} \|_{B_{p,r}^{s-1}}} + {\| {{\zeta^{(2)}}} \|_{B_{p,r}^{s-1}}}),
\end{eqnarray*}
and
\begin{eqnarray*}
&{\| {P(D)\big({\zeta ^{(12)}}(a_1({\zeta ^{(1)}} + {\zeta ^{(2)}}){u^{(1)}}
+ (2{a_1} + {a_2}){u^{(1)}}+{a_5}( {({u^{(2)}})^2}({\zeta ^{(1)}} + {\zeta ^{(2)}} + 1)))\big)\|_{B_{p,r}^{s - 1}}}}\\
&\le C{\| {{\zeta^{(12)}}} \|_{B_{p,r}^{s - 2}}}\big({\| {{u^{(1)}}} \|_{B_{p,r}^s}}+ {\| {{u^{(2)}}} \|^2_{B_{p,r}^s}}
+({\| {{\zeta^{(1)}}} \|_{B_{p,r}^{s-1}}} + {\| {{\zeta^{(2)}}} \|_{B_{p,r}^{s-1}}})({\| {{u^{(1)}}} \|_{B_{p,r}^s}} + {\| {{u^{(2)}}} \|^2_{B_{p,r}^s}}),
\end{eqnarray*}
which implies
\begin{eqnarray}\label{F}
{\| {{F(t)}} \|_{B_{p,r}^{s - 1}}} &\le& C({\| {{u^{(12)}}} \|_{B_{p,r}^{s - 1}}} + {\| {{\zeta ^{(12)}}} \|_{B_{p,r}^{s - 2}}}) \nonumber \\
&&\times \big( 1+\sum\limits_{k = 1}^3 ({{{\| {{u^{(1)}}} \|}_{B_{p,r}^s}} + {{\| {{u^{(2)}}} \|}_{B_{p,r}^s}} + {{\| {{\zeta ^{(1)}}} \|}_{B_{p,r}^{s - 1}}} + {{\| {{\zeta ^{(2)}}} \|}_{B_{p,r}^{s - 1}}}){^k}}\big).
\end{eqnarray}
Applying Lemma \ref{Lem6.2} $(i)$ to the first equation of the system (\ref{two-component}), we get
\begin{eqnarray}\label{u12}
{\| {{u^{(12)}}} \|_{B_{p,r}^{s - 1}}}&\le& {\| {u_0^{(12)}} \|_{B_{p,r}^{s - 1}}}e^ { C\int_0^t {{{\| {{u^{(2)}}(\tau )} \|}_{B_{p,r}^{s}}}} d\tau } \nonumber \\
&&+\int_0^t {e^ {C\int_\tau ^t {{{\| {{u^{(2)}}(\tau ')} \|}_{B_{p,r}^{s }}}} d\tau '}} \cdot {\| {{F}(\tau)} \|_{B_{p,r}^{s - 1}}}d\tau.
\end{eqnarray}
Similarly, since $B_{p,r}^{s-2} $ with $s>2 +\frac{1}{p}$ is an algebra, one can easily obtain
\begin{eqnarray}\label{G}
{\| {{G(t)}} \|_{B_{p,r}^{s - 2}}} &\le& C({\| {{u^{(12)}}} \|_{B_{p,r}^{s - 1}}} + {\| {{\zeta ^{(12)}}} \|_{B_{p,r}^{s - 2}}}) \nonumber \\
&&\times \big( 1+\sum\limits_{k = 1}^3 ({{{\| {{u^{(1)}}} \|}_{B_{p,r}^s}} + {{\| {{u^{(2)}}} \|}_{B_{p,r}^s}} + {{\| {{\zeta ^{(1)}}} \|}_{B_{p,r}^{s - 1}}} + {{\| {{\zeta ^{(2)}}} \|}_{B_{p,r}^{s - 1}}}){^k}}\big).
\end{eqnarray}
For $s>2 +\frac{1}{p}$, noticing that the embedding $B_{p,r}^{s-2}\hookrightarrow B_{p,r}^{\frac{1}{p}}\cap L^\infty$ holds true.
Then applying Lemma \ref{Lem6.2} $(i)$ again to second equation of the system (\ref{two-component}), we have
\begin{eqnarray}\label{zeta12}
{\| {{\zeta ^{(12)}}} \|_{B_{p,r}^{s - 2}}} &\le& {\| {\zeta _0^{(12)}} \|_{B_{p,r}^{s - 2}}}e^ {C\int_0^t ({{{\| {{u^{(2)}}(\tau )} \|}_{B_{p,r}^s}} + \| {{\zeta ^{(2)}}(\tau )} \|_{B_{p,r}^{s - 1}}^2 + {{\| {{\zeta ^{(2)}}(\tau )} \|}_{B_{p,r}^{s - 1}}})} d\tau }\nonumber\\
 &&+ \int_0^t {e^ {C\int_\tau ^t ({{{\| {{u^{(2)}}(\tau ')} \|}_{B_{p,r}^s}} + \| {{\zeta ^{(2)}}(\tau ')} \|_{B_{p,r}^{s - 1}}^2 + {{\| {{\zeta ^{(2)}}(\tau ')} \|}_{B_{p,r}^{s - 1}}}} )d\tau ' }} \cdot {\| {{G}}(\tau) \|_{B_{p,r}^{s - 2}}}d\tau.
\end{eqnarray}
Plugging (\ref{F}) into (\ref{u12}), and (\ref{G}) into (\ref{zeta12}), we derive from the resulting equations that
\begin{equation*}
\begin{aligned}
&{{\mathop{\rm e}\nolimits} ^{ - C\int_0^t ({{{\| {{u^{(2)}(\tau)}} \|}_{B_{p,r}^s}} +\| {{\zeta ^{(2)}(\tau)}} \|_{B_{p,r}^{s - 1}}^2 + {{\| {{\zeta ^{(2)}(\tau)}} \|}_{B_{p,r}^{s - 1}}}}) d\tau}}({\| {{u^{(12)}}} \|_{B_{p,r}^{s - 1}}} + {\| {{\zeta ^{(12)}}} \|_{B_{p,r}^{s - 2}}})\\
& \le( {\| {u_0^{(12)}} \|_{B_{p,r}^{s - 1}}} + {\| {\zeta _0^{(12)}} \|_{B_{p,r}^{s - 2}}})+\int_0^t {{{\mathop{\rm e}\nolimits} ^{ - C\int_0^\tau ( {{{\| {{u^{(2)}(\tau')}} \|}_{B_{p,r}^s}} + \| {{\zeta ^{(2)}(\tau')}} \|_{B_{p,r}^{s - 1}}^2 + {{\| {{\zeta ^{(2)}(\tau')}} \|}_{B_{p,r}^{s - 1}}})d\tau '} }}}\\
 &\times({\| {{u^{(12)}}} \|_{B_{p,r}^{s - 1}}} + {\| {{\zeta ^{(12)}}} \|_{B_{p,r}^{s - 2}}}) \cdot \big( 1+\sum\limits_{k = 1}^3 ({{{\| {{u^{(1)}}} \|}_{B_{p,r}^s}} + {{\| {{u^{(2)}}} \|}_{B_{p,r}^s}} + {{\| {{\zeta ^{(1)}}} \|}_{B_{p,r}^{s - 1}}} + {{\| {{\zeta ^{(2)}}} \|}_{B_{p,r}^{s - 1}}}){^k}}\big) d\tau.
\end{aligned}
\end{equation*}
Hence, using the Gronwall inequality, we reach (\ref{uniqueness}). This completes the proof of Lemma \ref{Lem3.1}.
\end{proof}

Next, we use classical Friedrichs regularization method to construct the approximation solutions to the system (\ref{two-component}). Before this,
for $T>0,s\in \R$ and $1\leq p\leq \infty$, we define
\begin{equation*}
E^s_{p,r}(T):=C([0,T];B^s_{p,r})\cap C^1([0,T];B^{s-1}_{p,r}), \quad
\mbox{if} \ r<\infty,
\end{equation*}
\begin{equation*}
E^s_{p,\infty}(T):=L^\infty(0,T;B^s_{p,\infty})\cap
\mbox{Lip}([0,T];B^{s-1}_{p,\infty}).
\end{equation*}
\begin{lemma3}\label{Lem3.2}
Let $1\leq p,r\leq \infty$ and $s>2 + \frac{1}{p}$. If $(u_0,\zeta_0)\in B_{p,r}^s\times B_{p,r}^{s-1}$, and $(u^{(0)},\zeta^{(0)}):=(0,0)$, then there exists a sequence of smooth functions $(u^{(j)},\zeta^{(j)})_{j\in \mathbb{N}}\in C(R^{+};B_{p,r}^\infty)^2$ solving the following linear transport equation by induction:
\begin{equation}\label{lineareq}
\left\{
\begin{aligned}
&{\partial _t}{u^{(j+1)}} + \sigma {u^{(j)}}{\partial _x}{u^{(j+1)}}=\tilde{F}(t,x), \\
&{\partial _t}{\zeta ^{(j+1)}} + ({u^{(j)}} + {b_1}{({\zeta ^{(j)}})^2} + 2{b_1}{\zeta ^{(j)}}){\partial _x}{\zeta ^{(j+1)}}=\tilde{G}(t,x), \\
&{u^{(j+1)}}{{\rm{|}}_{t = 0}} = u_0^{(j+1)}(x) = {S_{j+1}}{u_{0,}}                                                                                                                                                                                                     \\
&{\zeta ^{(j+1)}}{{\rm{|}}_{t = 0}} = \zeta _0^{(j+1)}(x) = {S_{j+1}}{\zeta _{0}},
\end{aligned}
\right.
\end{equation}
where $\tilde{F}(t,x):= P(D)\big(\frac{{3 - \sigma }}{2}{({u^{(j)}})^2} + \frac{\sigma }{2}{\partial _x}{({u^{(j)}})^2} + \frac{1}{2}{({\zeta ^{(j)}})^2} + {\zeta ^{(j)}} + {a_1}{u^{(j)}}({({\zeta ^{(j)}})^2} + 2{\zeta ^{(j)}} + 1)+ {a_2}{u^{(j)}}({\zeta ^{(j)}} + 1)+ ({a_3} + {a_4}){u^{(j)}} + {a_5}{({u^{(j)}})^2}{\zeta ^{(j)}}({\zeta ^{(j)}} + 1) + {a_6}{({u^{(j)}})^3}\big),$ and
$\tilde{G }(t,x):= - {\zeta ^{(j)}}{\partial _x}{u^{(j)}} - {\partial _x}{u^{(j)}} - {b_2}{\zeta ^{(j)}}{u^{(j)}}{\partial _x}{u^{(j)}} - {b_2}{u^{(j)}}{\partial _x}{u^{(j)}} - {b_3}{\zeta ^{(j)}}{({u^{(j)}})^2}{\partial _x}{u^{(j)}}- {b_3}{({u^{(j)}})^2}{\partial _x}{u^{(j)}}.$
Moreover, there exists some $T>0$ such that the solution $(u^{(j)},\zeta^{(j)})_{j\in \mathbb{N}}$ is uniformly bounded in $E_{p,r}^s(T)\times E_{p,r}^{s-1}(T)$,
and a Cauchy sequence in $C([0,T];B_{p,r}^{s-1})\times C([0,T];B_{p,r}^{s-2})$.
\end{lemma3}
\begin{proof}
Thanks to $S_{j+1}u_0, S_{j+1}\zeta_0 \in B_{p,r}^\infty$, it then follows Lemma \ref{Lem6.3} and by induction with respect to the index $j$ that the system (\ref{lineareq})
has a global solution which belongs to $C(R^{+};B_{p,r}^{\infty})^2$.

Using a similar argument in the proof of Lemma \ref{Lem3.1}, we have
\begin{equation}\label{estimate}
\begin{aligned}
&{\| {{u^{(j+1)}}} \|_{B_{p,r}^s}} + {\| {{\zeta ^{(j+1)}}} \|_{B_{p,r}^{s - 1}}}\le ({\| {{u_0}} \|_{B_{p,r}^s}} + {\| {{\zeta _0}} \|_{B_{p,r}^{s - 1}}})e^ {C\int_0^t \big({\sum\limits_{k = 1}^2 {{{({{\| {{u^{(j)}}(\tau )} \|}_{B_{p,r}^s}} + {{\| {{\zeta ^{(j)}}(\tau )} \|}_{B_{p,r}^{s - 1}}})}^k}} }\big) d\tau } \\
&+ \frac{C}{2}\int_0^t {e^ {C\int_\tau ^t \big({\sum\limits_{k = 1}^2 {{{({{\| {{u^{(j)}}({\tau '} )} \|}_{B_{p,r}^s}} + {{\| {{\zeta ^{(j)}}({\tau '} )} \|}_{B_{p,r}^{s - 1}}})}^k}} } \big)d\tau '}} \cdot \big(\sum\limits_{k = 1}^4 {{{({{\| {{u^{(j)}}(\tau )}\|}_{B_{p,r}^s}} + {{\| {{\zeta ^{(j)}}(\tau )} \|}_{B_{p,r}^{s - 1}}})}^k}} \big)d\tau.
\end{aligned}
\end{equation}
Choosing $$0<T \le \frac{1}{{3C\big(4({{\| {{u_0}} \|}_{B_{p,r}^s}} + {{\| {{\zeta _0}} \|}_{B_{p,r}^{s - 1}}}) + \sum\limits_{k = 1,k \ne 2}^4 {{2^{2k - 2}}{{({{\| {{u_0}} \|}_{B_{p,r}^s}} + {{\| {{\zeta _0}} \|}_{B_{p,r}^{s - 1}}})}^{k - 1}}} \big)}},$$
and supposing by induction that for all $t\in [0, T]$
\begin{equation}\label{induction}
{\| {{u^{(j)}}(t )} \|_{B_{p,r}^s}} + {\| {{\zeta ^{(j)}}(t )} \|_{B_{p,r}^{s - 1}}} \le \frac{{2({{\| {{u_0}} \|}_{B_{p,r}^s}} + {{\| {{\zeta _0}}\|}_{B_{p,r}^{s - 1}}})}}{{1 - 4C({{\| {{u_0}} \|}_{B_{p,r}^s}} + {{\| {{\zeta _0}} \|}_{B_{p,r}^{s - 1}}})t}}.
\end{equation}
Assume that (\ref{induction}) is valid for $j$, by using the mean-value theorem of integrals, we get for $\xi \in (\tau,t)$
\begin{eqnarray}\label{integral}
&&e^{ C\int_\tau ^t\big( {\sum\limits_{k = 1}^2 {{{({{\| {{u^{(j)}}({\tau '})} \|}_{B_{p,r}^s}} + {{\| {{\zeta ^{(j)}}({\tau '})} \|}_{B_{p,r}^{s - 1}}})}^k}} }\big) d\tau '}
\le e^{C\int_\tau ^t
\big( {\sum\limits_{k = 1}^2
{\frac{{2^k({{\| {{u_0}} \|}_{B_{p,r}^s}} + {{\| {{\zeta _0}} \|}_{B_{p,r}^{s - 1}}})^k}}
{({1 - 4C({{\| {{u_0}} \|}_{B_{p,r}^s}} + {{\| {{\zeta _0}} \|}_{B_{p,r}^{s - 1}}})\cdot{\tau '}})^k}
\big)}}d\tau '}\nonumber \\
 &=& e^ {\frac{1}{2}\ln\frac{{1 - 4C({{\| {{u_0}} \|}_{B_{p,r}^s}} + {{\| {{\zeta _0}} \|}_{B_{p,r}^{s - 1}}})\tau }}{{1 - 4C({{\| {{u_0}}\|}_{B_{p,r}^s}} + {{\| {{\zeta _0}} \|}_{B_{p,r}^{s - 1}}})t}} + \frac{{4C{{({{\| {{u_0}} \|}_{B_{p,r}^s}} + {{\| {{\zeta _0}} \|}_{B_{p,r}^{s - 1}}})}^2}(t - \tau )}}{{{{(1 - 4C({{\| {{u_0}} \|}_{B_{p,r}^s}} + {{\| {{\zeta _0}} \|}_{B_{p,r}^{s - 1}}})\xi )}^2}}}}\nonumber\\
 &\le& e^ {16C{({\| {{u_0}} \|_{B_{p,r}^s}} + {\| {{\zeta _0}} \|_{B_{p,r}^{s - 1}}})^2}T} \cdot{(\frac{{1 - 4C({{\| {{u_0}} \|}_{B_{p,r}^s}} + {{\| {{\zeta _0}} \|}_{B_{p,r}^{s - 1}}})\tau }}{{1 - 4C({{\| {{u_0}} \|}_{B_{p,r}^s}} + {{\| {{\zeta _0}} \|}_{B_{p,r}^{s - 1}}})t}})^{\frac{1}{2}}}\nonumber\\
& \le& {e^{\frac{1}{3}}}{(\frac{{1 - 4C({{\left\| {{u_0}} \right\|}_{B_{p,r}^s}} + {{\left\| {{\zeta _0}} \right\|}_{B_{p,r}^{s - 1}}})\tau }}{{1 - 4C({{\left\| {{u_0}} \right\|}_{B_{p,r}^s}} + {{\left\| {{\zeta _0}} \right\|}_{B_{p,r}^{s - 1}}})t}})^{\frac{1}{2}}}.
\end{eqnarray}
Then combining (\ref{induction}) and (\ref{integral}), by using the mean-value theorem of integrals again, we obtain
\begin{eqnarray}\label{estimate1}
&&\frac{C}{2}\int_0^t {e^ {C\int_\tau ^t \big({\sum\limits_{k = 1}^2 {{{({{\| {{u^{(j)}}({\tau '} )} \|}_{B_{p,r}^s}} + {{\| {{\zeta ^{(j)}}({\tau '} )} \|}_{B_{p,r}^{s - 1}}})}^k}} } \big)d\tau '}} \cdot \big(\sum\limits_{k = 1}^4 {{{({{\| {{u^{(j)}}(\tau )}\|}_{B_{p,r}^s}} + {{\| {{\zeta ^{(j)}}(\tau )} \|}_{B_{p,r}^{s - 1}}})}^k}} \big)d\tau\nonumber\\
 &\le& \frac{{{e^{\frac{1}{3}}}({{\| {{u_0}} \|}_{B_{p,r}^s}} + {{\| {{\zeta _0}} \|}_{B_{p,r}^{s - 1}}})}}{{{{(1 - 4C({{\| {{u_0}} \|}_{B_{p,r}^s}} + {{\| {{\zeta _0}} \|}_{B_{p,r}^{s - 1}}})t)}^{\frac{1}{2}}}}} \int_0^t {\big(\sum\limits_{k = 1}^4 {\frac{{{2^{k - 1}}C{{({{\| {{u_0}} \|}_{B_{p,r}^s}} + {{\| {{\zeta _0}} \|}_{B_{p,r}^{s - 1}}})}^{k - 1}}}}{{{{(1 - 4C({{\| {{u_0}} \|}_{B_{p,r}^s}} + {{\| {{\zeta _0}} \|}_{B_{p,r}^{s - 1}}})\tau )}^{k - \frac{1}{2}}}}}} } \big)d\tau \nonumber\\
& \le& \frac{{{e^{\frac{1}{3}}}({{\| {{u_0}} \|}_{B_{p,r}^s}} + {{\| {{\zeta _0}} \|}_{B_{p,r}^{s - 1}}})}}{{{{(1 - 4C({{\| {{u_0}} \|}_{B_{p,r}^s}} + {{\| {{\zeta _0}} \|}_{B_{p,r}^{s - 1}}})t)}^{\frac{1}{2}}}}} \cdot \big(\int_0^t {\frac{{2C({{\| {{u_0}} \|}_{B_{p,r}^s}} + {{\| {{\zeta _0}} \|}_{B_{p,r}^{s - 1}}})}}{{{{(1 - 4C({{\| {{u_0}} \|}_{B_{p,r}^s}} + {{\| {{\zeta _0}} \|}_{B_{p,r}^{s - 1}}})\tau )}^{\frac{3}{2}}}}}d\tau} \nonumber\\
&&+ \sum\limits_{k = 1,k \ne 2}^4 {\frac{{{2^{k - 1}}C{{({{\| {{u_0}} \|}_{B_{p,r}^s}} + {{\| {{\zeta _0}} \|}_{B_{p,r}^{s - 1}}})}^{k - 1}t}}}{{{{(1 - 4C({{\| {{u_0}} \|}_{B_{p,r}^s}} + {{\| {{\zeta _0}} \|}_{B_{p,r}^{s - 1}}})\xi_k )}^{k - \frac{1}{2}}}}}} \big)\nonumber\\
& \le& \frac{{{e^{\frac{1}{3}}}({{\| {{u_0}} \|}_{B_{p,r}^s}} + {{\| {{\zeta _0}} \|}_{B_{p,r}^{s - 1}}})}}{{{{(1 - 4C({{\| {{u_0}} \|}_{B_{p,r}^s}} + {{\| {{\zeta _0}} \|}_{B_{p,r}^{s - 1}}})t)}^{\frac{1}{2}}}}} \cdot \big((\frac{1}{{{{(1 - 4C({{\| {{u_0}} \|}_{B_{p,r}^s}} + {{\| {{\zeta _0}} \|}_{B_{p,r}^{s - 1}}})t )}^{\frac{1}{2}}}}} - 1)\nonumber\\
 &&+\sum\limits_{k = 1,k \ne 2}^4 {\frac{{{2^{k - 1}}C{{({{\| {{u_0}} \|}_{B_{p,r}^s}} + {{\| {{\zeta _0}} \|}_{B_{p,r}^{s - 1}}})}^{k - 1}}T}}{{{{(1 - 4C({{\| {{u_0}}. \|}_{B_{p,r}^s}} + {{\| {{\zeta _0}} \|}_{B_{p,r}^{s - 1}}}){t})}^{k- \frac{1}{2}}}}}} \big),
\end{eqnarray}
where $0<\xi_{k}<t$, $(1\leq k\leq 4,k\neq 2)$. Inserting (\ref{integral}) with $\tau =0 $ and (\ref{estimate1}) into (\ref{estimate}), we get
\begin{eqnarray}\label{estimatej+1}
&&{\| {{u^{(j+1)}}} \|_{B_{p,r}^s}} + {\| {{\zeta ^{(j+1)}}}\|_{B_{p,r}^{s - 1}}} \nonumber\\
  &\le& \frac{{{e^{\frac{1}{3}}}({{\| {{u_0}} \|}_{B_{p,r}^s}} + {{\| {{\zeta _0}} \|}_{B_{p,r}^{s - 1}}})}}{{(1 - 4C({{\| {{u_0}} \|}_{B_{p,r}^s}} + {{\| {{\zeta _0}} \|}_{B_{p,r}^{s - 1}}})t)}}\big(1 + \sum\limits_{k= 1,k \ne 2}^4 {\frac{{{2^{k - 1}}C{{({{\| {{u_0}} \|}_{B_{p,r}^s}} + {{\| {{\zeta _0}} \|}_{B_{p,r}^{s - 1}}})}^{k- 1}}T}}{{{{(1 - 4C({{\| {{u_0}} \|}_{B_{p,r}^s}} + {{\| {{\zeta _0}} \|}_{B_{p,r}^{s - 1}}}){t})}^{k - 1}}}}} \big)\nonumber\\
 &\le& \frac{{{e^{\frac{1}{3}}}({{\| {{u_0}} \|}_{B_{p,r}^s}} + {{\| {{\zeta _0}} \|}_{B_{p,r}^{s - 1}}})}}{{(1 - 4C({{\| {{u_0}} \|}_{B_{p,r}^s}} + {{\| {{\zeta _0}} \|}_{B_{p,r}^{s - 1}}})t)}}(1 + CT\sum\limits_{k = 1,k \ne 2}^4 {{2^{2k - 2}}{{({{\| {{u_0}} \|}_{B_{p,r}^s}} + {{\| {{\zeta _0}} \|}_{B_{p,r}^{s - 1}}})}^{k - 1}}} )\nonumber\\
 &\le& \frac{{4{e^{\frac{1}{3}}}({{\left\| {{u_0}} \right\|}_{B_{p,r}^s}} + {{\left\| {{\zeta _0}} \right\|}_{B_{p,r}^{s - 1}}})}}{{3(1 - 4C({{\left\| {{u_0}} \right\|}_{B_{p,r}^s}} + {{\left\| {{\zeta _0}} \right\|}_{B_{p,r}^{s - 1}}})t)}} \le \frac{{2({{\left\| {{u_0}} \right\|}_{B_{p,r}^s}} + {{\left\| {{\zeta _0}} \right\|}_{B_{p,r}^{s - 1}}})}}{{1 - 4C({{\left\| {{u_0}} \right\|}_{B_{p,r}^s}} + {{\left\| {{\zeta _0}} \right\|}_{B_{p,r}^{s - 1}}})t}},
\end{eqnarray}
which implies that $(u^{(j)},\zeta^{(j)})_{j\in \mathbb{N}}$ is uniformly bounded in $C([0,T];B_{p,r}^s)\times C([0,T];B_{p,r}^{s-1})$.
Since $B_{p,r}^{s-1}, B_{p,r}^{s-2}$ with $s>2 + \frac{1}{p}$ is an algebra, one can readily prove
that $\sigma u^{(j)}\partial_x u^{(j+1)},\tilde{F}$ are uniformly bounded in $C([0,T];B_{p,r}^{s-1})$, and $(u^{(j)}+b_1(\zeta^{(j)})^2+2b_1\zeta^{(j)})\partial_x\zeta^{(j+1)}, \tilde{G}$ are uniformly bounded in $C([0,T];B_{p,r}^{s-2})$. It thus follows from the system (\ref{lineareq}) that $(\partial_t u^{(j+1)},\partial_t\zeta^{(j+1)})\in C([0,T];B_{p,r}^{s-1})\times C([0,T];B_{p,r}^{s-2})$.
This yields the uniform bound of $(u^{(j)},\zeta^{(j)})_{j\in \mathbb{N}}$ in $E_{p,r}^s(T)\times E_{p,r}^{s-1}(T)$.

Next we prove that $(u^{(j)},\zeta^{(j)})_{j\in \mathbb{N}}$ is a Cauchy sequence in $C([0,T];B_{p,r}^{s-1})\times C([0,T];B_{p,r}^{s-2})$. Indeed, for $j,l\in \mathbb{N}$, we deduce from
the system (\ref{lineareq}) that
\begin{equation}
\left\{
\begin{aligned}
&({\partial _t} + \sigma {u^{(j+l)}}{\partial _x})({u^{(j+l+1)}} - {u^{(j+1)}})= \bar{F}(t,x),\\
&({\partial _t} + {(u^{(j+l)}+b_1(\zeta^{(j+l)})^2+2b_1\zeta^{(j+l)})}{\partial _x})({\zeta ^{(j+l+1)}} - {\zeta ^{(j+1)}}) =\bar{G}(t,x),
\end{aligned}
\right.
\end{equation}
where $\bar{F}(t,x):=\sigma ({u^{(j)}} - {u^{(j+l)}}){\partial _x}{u^{(j+1)}}+P(D)\big(
({u^{(j+l)}} - {u^{(j)}})(\frac{{3 - \sigma }}{2}({u^{(j+l)}} + {u^{(j)}})+a_1{({\zeta ^{(j)}})^2}
+ (2{a_1}+a_2){\zeta ^{(j)}}+{(a_1+a_2+a_3+a_4)}+{a_5}({u^{(j+l)}} + {u^{(j)}})({({\zeta ^{(j)}})^2}+{{\zeta ^{(j)}}})
+{a_6}({({u^{(j+l)}})^2} + {u^{(j+l)}}{u^{(j)}} + {({u^{(j)}})^2}))
+ \frac{\sigma }{2}({\partial _x}{u^{(j+l)}} - {\partial _x}{u^{(j)}})({\partial _x}{u^{(j+l)}} + {\partial _x}{u^{(j)}})
+({\zeta ^{(j+l)}} - {\zeta ^{(j)}})((\frac{1}{2}+{a_1}{u^{(j+l)}})({\zeta ^{(j+l)}} + {\zeta ^{(j)}}+2)
+a_2{u^{(j+l)}}+{a_5}{({u^{(j+l)}})^2}({\zeta ^{(j+l)}} + {\zeta ^{(j)}}+1))
\big)$, and $\bar{G}(t,x):={\partial _x}({u^{(j)}} - {u^{(j+l)}})\big({\zeta ^{(j)}}+1+({b_2}{u^{(j)}}+{b_3}{({u^{(j)}})^2})({\zeta ^{(j+l)}}+1)\big)+({\zeta ^{(j)}} - {\zeta ^{(j+l)}})\big(
{\partial _x}{u^{(j+l)}}+{u^{(j)}}{\partial _x}{u^{(j)}}({b_2}+{b_3}{u^{(j)}})+
{b_1}{\partial _x}{\zeta ^{(j+1)}} ({\zeta ^{(j)}} + {\zeta ^{(j+l)}}+2)\big)
+({u^{(j)}} - {u^{(j+l)}})\big(({\zeta ^{(j+l)}}+1)
({b_2}{\partial _x}{u^{(j+l)}}+{b_3}({u^{(j)}} + {u^{(j+l)}}){\partial _x}{u^{(j+l)}})
+{\zeta ^{(j+1)}})
\big).$ Similarly, it then follows Lemma \ref{Lem6.2} $(i)$ and the fact that $B_{p,r}^{s-1},B_{p,r}^{s-2}\ (s>2+\frac{1}{p})$ are algebras that
\begin{eqnarray}\label{estimate2}
&&{\| {{u^{(j+l+1)}} - {u^{(j+1)}}} \|_{B_{p,r}^{s - 1}}}+{\| {{\zeta ^{(j+l+1)}} - {\zeta ^{(j+1)}}}\|_{B_{p,r}^{s - 2}}}\nonumber\\
&\le& {e^{  C{{\int_0^t {{{\| {{u^{(j+l)}}} \|}_{B_{p,r}^s}} + \| {{\zeta ^{(j+l)}}} \|_{B_{p,r}^{s - 1}}^2 +\| {{\zeta ^{(j+l)}}} \|} }_{B_{p,r}^{s - 1}}}d\tau }}
\big
(({\| {{u_0^{(j+l+1)}} - {u_0^{(j+1)}}} \|_{B_{p,r}^{s - 1}}}\nonumber\\
&&
+{\| {{\zeta_0 ^{(j+l+1)}} - {\zeta_0 ^{(j+1)}}}\|_{B_{p,r}^{s - 2}}})+\int_0^t{e^{ - C{{\int_0^\tau {{{\| {{u^{(j+1)}}} \|}_{B_{p,r}^s}} + \| {{\zeta ^{(j+1)}}} \|_{B_{p,r}^{s - 1}}^2 +\| {{\zeta ^{(j+1)}}} \|} }_{B_{p,r}^{s - 1}}}d\tau' }}\nonumber\\
  &&\times({\| {{u^{(j+l)}} - {u^{(j)}}} \|_{B_{p,r}^{s - 1}}}+{\| {{\zeta ^{(j+l)}} - {\zeta ^{(j)}}}\|_{B_{p,r}^{s - 2}}})\cdot\big(
  \sum\limits_{k = 1}^3({{\| {{u^{(j)}}} \|}_{B_{p,r}^s}}+{{\| {{u^{(j+1)}}} \|}_{B_{p,r}^s}} + {{\| {{u^{(j+l)}}} \|}_{B_{p,r}^s}}
  \nonumber\\
  &&\times { {{\| {{\zeta ^{(j)}}}\|}_{B_{p,r}^{s - 1}}} + {{\| {{\zeta ^{(j+1)}}} \|}_{B_{p,r}^{s-1}}} + {{\| {{\zeta ^{(j+l)}}} \|}_{B_{p,r}^{s - 1}}} + 1){^k}}\big) d\tau\big).
\end{eqnarray}
According to the fact that $(u^{(j)},\zeta^{(j)})_{j\in \mathbb{N}}$ is uniformly bounded in $E_{p,r}^s(T)\times E_{p,r}^{s-1}(T)$, and
${u_0^{(j+l+1)}} - {u_0^{(j+1)}}=\sum_{q = j+ 1}^{j+l} {{\Delta _q}{u_0}}$, ${\zeta_0^{(j+l+1)}} - {\zeta_0^{(j+1)}}=\sum_{q = j+ 1}^{j+l} {{\Delta _q}{\zeta_0}}$,
there exists a positive constant $C_T$ independent of $j, l$  such that
\begin{eqnarray*}
&&{\| {{u^{(j+l+1)}} - {u^{(j+1)}}} \|_{B_{p,r}^{s - 1}}}+{\| {{\zeta ^{(j+l+1)}} - {\zeta ^{(j+1)}}}\|_{B_{p,r}^{s - 2}}}\\
 &\le& {C_T}\big({2^{ - j}} + \int_0^t {({\| ({{u^{(j+l)}} - {u^{(j)}}})(\tau) \|_{B_{p,r}^{s - 1}}}+{\| ({{\zeta ^{(j+l)}} - {\zeta ^{(j)}}})(\tau)\|_{B_{p,r}^{s - 2}}})d\tau }\big).
\end{eqnarray*}
Arguing by induction with respect to the index $j$, one obtains
\begin{eqnarray*}
&&{\| {{u^{(j+l+1)}} - {u^{(j+1)}}} \|_{B_{p,r}^{s - 1}}}+{\| {{\zeta ^{(j+l+1)}} - {\zeta ^{(j+1)}}}\|_{B_{p,r}^{s - 2}}}\\
 &\le& {C_T}\big(2^{-j}\sum\limits_{k = 0}^j\frac{(2TC_T)^k}{k!}+C_T^{j+1}\int_0^t\frac{(t-\tau)^j}{j!}d\tau
\big)\\
&\le&{C_T}2^{-j}\sum\limits_{k = 0}^j\frac{(2TC_T)^k}{k!}+C_T\frac{(TC_T)^{j+1}}{(j+1)!},
\end{eqnarray*}
which yields the desired result as $j\rightarrow \infty$. This completes the proof of Lemma \ref{Lem3.2}.
\end{proof}

Now, we are going to complete the proof of the main theorem of this section.
\begin{theorem3}\label{Th3.1}
Assume that $1 \leq p,r \leq \infty$, $s > 2+ \frac{1}{p}$, and the initial data $(u_0,\zeta_0)\in B_{p,r}^s\times B_{p,r}^{s-1}$. Then there exists a time $T>0$ such that the Cauchy problem of the system $(\ref{two-component})$
has a unique solution $(u, \zeta) \in E_{p,r}^s(T)\times E_{p,r}^{s-1}(T)$, and the mapping $(u_0, \zeta_0)\mapsto(u, \zeta)$ is continuous from a neighborhood of $B_{p,r}^s\times B_{p,r}^{s-1}$ into
$C([0,T];B_{p,r}^{s^\prime})\cap C^{1}([0,T];B_{p,r}^{s^\prime-1})\times C([0,T];B_{p,r}^{s^\prime-1})\cap C^{1}([0,T];B_{p,r}^{s^\prime-2})$ for all $s'<s$ when $r=\infty$, and $s'=s$ whereas $r<\infty$.
\end{theorem3}
\begin{proof}
According to Lemma \ref{Lem3.2}, we denote $(u,\zeta)$ as the limit of the existing Cauchy sequence $(u^{(j)},\zeta^{(j)})_{j\in \mathbb{N}}$ in $C([0,T];B_{p,r}^{s-1})\times C([0,T];B_{p,r}^{s-2})$.
Now we check that $(u,\zeta)\in E_{p,r}^{s}(T)\times E_{p,r}^{s-1}(T)$ solves the system (\ref{two-component}). Indeed, combining the uniform boundedness of $(u^{(j)},\zeta^{(j)})_{j\in \mathbb{N}}$ in
$L^{\infty}(0,T;\\B_{p,r}^{s})\times L^{\infty}(0,T;B_{p,r}^{s-1})$ and the Fatou lemma (Lemma \ref{Lem6.1} $(ii)$), we have $(u,\zeta)\in L^{\infty}(0,T;B_{p,r}^{s})\times L^{\infty}(0,T;B_{p,r}^{s-1})$.
Then by the interpolation inequality (Lemma \ref{Lem6.1} $(iii)$), we get $(u^{(j)},\zeta^{(j)})_{j\in \mathbb{N}}\\ \rightarrow (u,\zeta)$, as $j\rightarrow\infty$, in $C([0,T];B_{p,r}^{s'})\times C([0,T];B_{p,r}^{s'-1})$, for all $s'<s$. Taking limit in (\ref{lineareq}), we conclude that $(u,\zeta)$ solves the system (\ref{two-component}) in the sense of $C([0,T];B_{p,r}^{s'-1})\times C([0,T];B_{p,r}^{s'-2})$, for all $s'<s$.
Taking account of the system (\ref{two-component}) and Lemma \ref{Lem6.3}, for $r<+\infty$, we infer that  $(u,\zeta)\in C([0,T];B_{p,r}^{s'})\times C([0,T];B_{p,r}^{s'-1})$ for all $s'\leq s$. Using
the system (\ref{two-component}) again, we see that $(\partial_tu,\partial_t\zeta)\in C([0,T];B_{p,r}^{s-1})\times C([0,T];B_{p,r}^{s-2})$ for $r<+\infty$, and in $L^{\infty}(0,T;B_{p,r}^{s-1})\times L^{\infty}(0,T;B_{p,r}^{s-2})$ otherwise. Hence $(u,\zeta)\in E_{p,r}^{s}(T)\times E_{p,r}^{s-1}(T)$. Moreover, the continuity with respect to the initial data in
$C([0,T];B_{p,r}^{s'})\cap C^{1}([0,T];B_{p,r}^{s'-1})\times C([0,T];B_{p,r}^{s'-1})\cap C^{1}([0,T];B_{p,r}^{s'-2})$ for all $s'<s$ can be proved by Lemma \ref{Lem3.1} and an interpolation argument.
While the continuity up to $s'=s$ for $r<\infty$ can be obtained by making use of a sequence of viscosity approximation of solutions $(u_\varepsilon,\zeta_\varepsilon)_{\varepsilon>0}$ to the system (\ref{two-component}), which converges uniformly in $C([0,T];B_{p,r}^{s})\cap C^{1}([0,T];B_{p,r}^{s-1})\times C([0,T];B_{p,r}^{s-1})\cap C^{1}([0,T];B_{p,r}^{s-2})$. This completes the proof of Theorem \ref{Th3.1}.
\end{proof}

\begin{remark3}\label{Rem3.1}
It is noted to point that $\zeta=\rho-1$ in $(\ref{two-component0})$ is transported along direction of $u$ instead of the complicated $(u+b_1\zeta^2+2b_1\zeta)$ in $(\ref{two-component})$.
Going along the lines of the proof of Theorem $\ref{Th3.1}$,
one may follow the similar argument as in $\cite{G-L}$ to obtain
the local well-posedness result for the system $(\ref{two-component0})$ in $B_{p,r}^s\times B_{p,r}^{s-1}$ with $1 \leq p,r \leq \infty$ and $s >\max\{ 1+ \frac{1}{p},\frac{3}{2}\}.$
As is well-known, when $p=r=2$, $B_{p,r}^s(\R)=H^s({\R})$. Theorem $\ref{Th3.1}$ gives the local well-posedness for the system $(\ref{two-component})$ with initial data
$(u_0,\rho_0-1)\in H^s(\R)\times H^{s-1}(\R)$, $s>\frac{5}{2}$. While for the system $(\ref{two-component0})$, according to the above comments, it is locally well-posed in
$H^s(\R)\times H^{s-1}(\R)$, $s>\frac{3}{2}.$ Thus we need here higher regularity index to overcome the difficulties encountered by the second equation in $(\ref{two-component})$.
\end{remark3}

\section{Blow-up criterion}\label{sec4}
\newtheorem{theorem4}{Theorem}[section]
\newtheorem{lemma4}{Lemma}[section]
\newtheorem {remark4}{Remark}[section]
\newtheorem {definition4}{Definition}[section]
\newtheorem{corollary4}{Corollary}[section]
\par

In this section, we present the following blow-up criterion for the the general system (\ref{two-component}) in nonhomogeneous Besov spaces.
\begin{theorem4}\label{Th4.1}
Given $(u_0,\zeta_0)\in B_{p,r}^s\times B_{p,r}^{s-1}$ with $1 \leq p,r \leq \infty$ and $s > 2+ \frac{1}{p}$. Let $T$ be the maximal existence time of the corresponding solution $(u, \zeta)$ to the system $(\ref{two-component})$ with $(u_0,\zeta_0)\in B_{p,r}^s\times B_{p,r}^{s-1}$. If $T<\infty$, then
$$\int_0^{T}\big(\sum\limits_{k = 1}^2(\|u\|^k_{L^\infty}+\|\zeta\|^k_{L^\infty})+\|u_x\|_{L^\infty}+\|\zeta_x\|_{L^\infty}(1+\|\zeta\|_{L^\infty})
+\|u\|^2_{L^\infty}\|\zeta\|_{L^\infty}\big)(\tau)d\tau=\infty.$$
\end{theorem4}
\begin{proof}
Applying $\Delta_q$ to the system (\ref{two-component}), we have
\begin{equation}\label{4.1}
\left\{
 \begin{aligned}
&\Delta_qu_t+\sigma \Delta_qu_xu=\sigma[u,\Delta_q]\partial_xu+\Delta_q \hat{F},\\
&\Delta_q\zeta_t+\Delta_q\zeta_x (u+b_1\zeta^2+2b_1\zeta)=[u+b_1\zeta^2+2b_1\zeta,\Delta_q]\partial_x\zeta+\Delta_q \hat{G},
\end{aligned}
\right.
\end{equation}
where $[\cdot,\cdot]$ denotes the commutator of the operators, $\hat{F}:=P(D)\big(\frac{3-\sigma}{2}u^2+\frac{\sigma}{2}(u_x)^2+\frac{1}{2}\zeta^2+\zeta+a_1(\zeta^2 u+2\zeta u+u)
+a_2(\zeta u+u)+(a_3+a_4)u+a_5(u^2\zeta^2+u^2\zeta)+a_6u^3\big)$, and $\hat{G}:=-\zeta u_x-u_x-b_2\zeta uu_x-b_2uu_x-b_3\zeta u^2u_x-b_3u^2u_x$. Multiplying both sides of the first equation of (\ref{4.1}) by $\mbox{sgn}(\Delta_q u)|\Delta_q u|^{p-1}$ and integrating over $\R$ with respect to $x$, we obtain
\begin{eqnarray*}
\frac{1}{p}\frac{d}{dt}\int_{\R}|\Delta_qu|^pdx-\frac{\sigma}{p}\int_{\R}|\Delta_qu|^pu_xdx=\int_{\R}\mbox{sgn}(\Delta_q u)|\Delta_q u|^{p-1}(\sigma[u,\Delta_q]\partial_xu+\Delta_q \hat{F})dx.
\end{eqnarray*}
It then follows from H\"{o}lder's inequality and integration with respect to $t$ that
\begin{eqnarray}\label{4.2}
\|\Delta_qu\|_{L^p}\leq \|\Delta_qu_0\|_{L^p}+C\int_0^t(\|\Delta_qu\|_{L^p}\|u_x\|_{L^\infty}+\|[u,\Delta_q]\partial_xu\|_{L^p}+\|\Delta_q \hat{F}\|_{L^p})(\tau)d\tau.
\end{eqnarray}
Applying the similar method to the second equation of (\ref{4.1}), we deduce that
\begin{eqnarray}\label{4.3}
\|\Delta_q\zeta\|_{L^p}&\leq& \|\Delta_q\zeta_0\|_{L^p}+C\int_0^t\big(\|\Delta_q\zeta\|_{L^p}(\|u_x\|_{L^\infty}+\|\zeta\|_{L^\infty}\|\zeta_x\|_{L^\infty}+\|\zeta_x\|_{L^\infty})\nonumber\\
&&+\|[u+b_1\zeta^2+2b_1\zeta,\Delta_q]\partial_x\zeta\|_{L^p}+\|\Delta_q \hat{G}\|_{L^p}\big)(\tau)d\tau.
\end{eqnarray}
Multiplying both sides of (\ref{4.2}) by $2^{qs}$, taking $l^r$-norm and using Minkowski's inequality, we get
\begin{eqnarray}\label{4.4}
\|u\|_{B_{p,r}^s}\leq \|u_0\|_{B_{p,r}^s}+C\int_0^t(\|u\|_{B_{p,r}^s}\|u_x\|_{L^\infty}+\|2^{qs}\|[u,\Delta_q]\partial_xu\|_{L^p}\|_{l^r}+\|\hat{F}\|_{B_{p,r}^s})(\tau)d\tau.
\end{eqnarray}
Making use of Lemma \ref{Lem6.4}, we have
\begin{eqnarray}\label{4.5}
\|2^{qs}\|[u,\Delta_q]\partial_xu\|_{L^p}\|_{l^r}\leq C \|u\|_{B_{p,r}^s}\|u_x\|_{L^\infty}.
\end{eqnarray}
Thanks to the fact that $P(D)\in Op(S^{-1})$ and $B_{p,r}^{s-1}$ is an algebra, we infer from Lemma \ref{Lem6.1} that
\begin{eqnarray}\label{4.6}
\|\hat{F}\|_{B_{p,r}^s} \leq C \big(1+\sum\limits_{k = 1}^2(\|u\|^k_{L^\infty}+\|\zeta\|^k_{L^\infty})+\|u_x\|_{L^\infty}+\|u\|^2_{L^\infty}\|\zeta\|_{L^\infty}\big)(\|u\|_{B_{p,r}^s}+\|\zeta\|_{B_{p,r}^{s-1}}).
\end{eqnarray}
Similarly, multiplying both sides of (\ref{4.3}) by $2^{q(s-1)}$, taking $l^r$-norm and using Minkowski's inequality, it yields
\begin{eqnarray}\label{4.7}
\|\zeta\|_{B_{p,r}^{s-1}}&\leq& \|\zeta_0\|_{B_{p,r}^{s-1}}+C\int_0^t\big(\|\zeta\|_{B_{p,r}^{s-1}}(\|u_x\|_{L^\infty}+\|\zeta\|_{L^\infty}\|\zeta_x\|_{L^\infty}+\|\zeta_x\|_{L^\infty})\nonumber\\
&&+\|2^{q(s-1)}\|[u+b_1\zeta^2+2b_1\zeta,\Delta_q]\partial_x\zeta\|_{L^p}\|_{l^r}+\|\hat{G}\|_{B_{p,r}^{s-1}}\big)(\tau)d\tau.
\end{eqnarray}
By virtue of Lemma \ref{Lem6.4}, we obtain
\begin{eqnarray}\label{4.8}
&&\|2^{q(s-1)}\|[u+b_1\zeta^2+2b_1\zeta,\Delta_q]\partial_x\zeta\|_{L^p}\|_{l^r}\nonumber\\
&\leq& C \big(\|\zeta\|_{B_{p,r}^{s-1}}(\|u_x\|_{L^\infty}+\|\zeta\|_{L^\infty}\|\zeta_x\|_{L^\infty}+\|\zeta_x\|_{L^\infty})+\|u\|_{B_{p,r}^s}\|\zeta_x\|_{L^\infty}\big).
\end{eqnarray}
Noticing that $B_{p,r}^{s-1}$ is an algebra, we get
\begin{eqnarray}\label{4.9}
\|\hat{G}\|_{B_{p,r}^{s-1}} \leq C (1+\|u\|_{L^\infty}+\|u\|^2_{L^\infty}+\|\zeta\|_{L^\infty}+\|u\|_{L^\infty}\|\zeta\|_{L^\infty}+\|u\|^2_{L^\infty}\|\zeta\|_{L^\infty}
)\|u\|_{B_{p,r}^s}.
\end{eqnarray}
Now plugging (\ref{4.5})-(\ref{4.6}) into (\ref{4.4}) and (\ref{4.8})-(\ref{4.9}) into (\ref{4.7}) respectively, we have
\begin{eqnarray*}
&&\|u\|_{B_{p,r}^s}+\|\zeta\|_{B_{p,r}^{s-1}}\leq \|u_0\|_{B_{p,r}^s}+\|\zeta_0\|_{B_{p,r}^{s-1}}+C\int_0^t(\|u\|_{B_{p,r}^s}+\|\zeta\|_{B_{p,r}^{s-1}})\nonumber\\
&&\times\big(1+\sum\limits_{k = 1}^2(\|u\|^k_{L^\infty}+\|\zeta\|^k_{L^\infty})+\|u_x\|_{L^\infty}+\|\zeta_x\|_{L^\infty}(1+\|\zeta\|_{L^\infty})
+\|u\|^2_{L^\infty}\|\zeta\|_{L^\infty}\big)(\tau)d\tau.
\end{eqnarray*}
Taking advantage of Gronwall's inequality, we obtain
\begin{eqnarray}\label{4.10}
&&\|u\|_{B_{p,r}^s}+\|\zeta\|_{B_{p,r}^{s-1}}
\leq (\|u_0\|_{B_{p,r}^s}+\|\zeta_0\|_{B_{p,r}^{s-1}})\nonumber\\
&&\times e^{C\int_0^t\big(1+\sum\limits_{k = 1}^2(\|u\|^k_{L^\infty}+\|\zeta\|^k_{L^\infty})+\|u_x\|_{L^\infty}+\|\zeta_x\|_{L^\infty}(1+\|\zeta\|_{L^\infty})
+\|u\|^2_{L^\infty}\|\zeta\|_{L^\infty}\big)(\tau)d\tau}.
\end{eqnarray}
Hence, if $T<\infty$ satisfies that
$\int_0^{T}\big(\sum\limits_{k = 1}^2(\|u\|^k_{L^\infty}+\|\zeta\|^k_{L^\infty})+\|u_x\|_{L^\infty}+\|\zeta_x\|_{L^\infty}(1+\|\zeta\|_{L^\infty})
+\|u\|^2_{L^\infty}\|\zeta\|_{L^\infty}\big)(\tau)d\tau<\infty$, then we deduce from (\ref{4.10}) that
\begin{eqnarray*}
\limsup\limits_{t\rightarrow T}(\|u\|_{B_{p,r}^s}+\|\zeta\|_{B_{p,r}^{s-1}})<\infty,
\end{eqnarray*}
which contradicts the fact that $T$ is the maximal existence time. This completes the proof of Theorem \ref{Th4.1}.
\end{proof}

\begin{remark4}\label{Rem4.1}
When $A=0$, on the one hand, the conserved quantity $E(u(t),\zeta(t)):=\int_{\R}(u^2+u_x^2+\zeta^2)dx$ is independent of $t$ $\cite{H-L}$, which gives the uniform boundedness of $u$ by the Sobolev embedding.
On the other hand, the variable $\zeta=\rho-1$ is advected by the fluid flow $u$ in $(\ref{two-component0})$, which is different from the flow generated by $(u+b_1\zeta^2+2b_1\zeta)$ in the system $(\ref{two-component})$.
Thus we do not need the estimates related to $\|u\|_{L^\infty}$, $\|\zeta\|_{L^\infty}$ and $\|\zeta_x\|_{L^\infty}$ to obtain the blow-up criterion
for $(\ref{two-component0})$. Therefore, following the same procedure in the proof of Theorem $\ref{Th4.1}$, for $(u_0,\zeta_0)\in B_{p,r}^s\times B_{p,r}^{s-1}$ with $1 \leq p,r \leq \infty$ and $s >\max\{ 1+ \frac{1}{p},\frac{3}{2}\},$ if the maximal time of existence $T<\infty$, then $\int_0^{T}\|u_x(\tau)\|_{L^\infty}d\tau=\infty.$ When $p=r=2$, it coincides with the blow-up result in $\cite{Dong}$.
\end{remark4}

\section{Global existence when $A=0$ and $\sigma=0$}\label{sec5}
\newtheorem{theorem5}{Theorem}[section]
\newtheorem{lemma5}{Lemma}[section]
\newtheorem {remark5}{Remark}[section]
\newtheorem {definition5}{Definition}[section]
\newtheorem{corollary5}{Corollary}[section]
\par
Note that the wave-breaking phenomena for the system (\ref{two-component0}) (i.e., the vorticity $A$ vanishes in (\ref{2CHvor}))
when $\sigma\neq 0$ was established in \cite{Dong}. In this section, motivated by \cite{CLQ}, we present a sufficient condition for the global solution of the system (\ref{two-component0})
when $\sigma=0.$ To this end, we first rewrite (\ref{two-component0}) with $\sigma=0$ as the following form
\begin{equation}\label{two-componentnew}
\left\{
 \begin{aligned}
&{u_t} =-\partial_x p \ast ( \frac{3}{2}u^2 +\frac{1}{2}\rho^2-3u^2(\rho^2-\rho)+4u^3),\\
&{\rho _t} + {(\rho u)_x}=\rho {({u^3})_x},
\end{aligned}
\right.
\end{equation}
where $p(x):=\frac{1}{2}e^{-|x|}$, $x\in\R$.

To use the method of characteristics, we consider the following associated Lagrangian flow
\begin{equation}\label{flow}
\begin{aligned}
\left\{ {\begin{array}{*{20}{l}}
q_t(t,x) =  u(t,q(t,x)),   &t \in [0,T),\\
q(0,x) = x,  & x \in \R,
\end{array}} \right.
\end{aligned}
\end{equation}
where $u\in C^1([0,T);H^{s-1})$ is the first component of the solution $(u,\rho)$ to (\ref{two-componentnew}) with initial data $(u_0,\rho_0)\in H^s\times H^{s-1}$ with $s>\frac{3}{2}$, and
$T>0$ is the maximal time of existence. It then follows from a direct calculation that $q_x(t,x)=\exp\{\int_0^tu_x(\tau,q(\tau,x))d\tau\}>0,$ for all $(t,x)\in [0,T)\times \R.$ Thus the mapping $q(t,\cdot):\R\mapsto \R$
is an increasing diffeomorphism of $\R$. Hence the $L^\infty$-norm of any function $v(t,\cdot)\in L^\infty$ is preserved under $q(t,\cdot)$, i.e., $\|v(t,\cdot)\|_{L^\infty}=\|v(t,q(t,\cdot))\|_{L^\infty}$, $t\in [0,T).$
Moreover, we have $\inf_{x\in \R} v(t,x)=\inf_{x\in \R} v(t,q(t,x))$ and $\sup_{x\in \R} v(t,x)=\sup_{x\in \R} v(t,q(t,x))$, $t\in [0,T).$ On the other hand, we also need the following useful result given in \cite{Constantin-E}.
\begin{lemma5}\label{Lem5.1}
Let $T>0$ and $v\in C^1([0,T);H^2(\R))$. Then for every $t\in [0,T)$ there exists at least one point $\xi(t)\in \R$ $($respectively, $\gamma(t)\in \R$$)$
with $m(t):=\inf\limits_{x\in \R}\{ v_x(t,x)\}=v_x(t,\xi(t))$$($respectively, $M(t):=\sup\limits_{x\in \R}\{v_x(t,x)\}=v_x(t,\gamma(t))$$)$
, and the function $m$ $($respectively, $M(t)$$)$
is almost everywhere differentiable with $\frac{dm(t)}{dt}=v_{tx}(t,\xi(t))$ $($respectively, $\frac{dM(t)}{dt}=v_{tx}(t,\gamma(t))$$)$ a.e. on $(0,T)$.
\end{lemma5}

Indeed, in view of the blow-up criterion for $(\ref{two-component0})$ discussed in Remark $\ref{Rem4.1}$,
to obtain global well-posedness of solutions for (\ref{two-componentnew}), the following estimates of $u_x$ is essential.
\begin{lemma5}\label{Lem5.2}
Assume $E(0):=E(u(0),\rho(0))<\frac{1}{3}$ $($$E(u(t),\rho(t))$ is the conservation law in Remark $\ref{Rem4.1}$$)$. Let $(u,\rho)$ be the solution of the system $(\ref{two-componentnew})$ with initial data $(u_0,\rho_0-1)\in H^s\times H^{s-1}$, $s>\frac{3}{2}$, and $T$ be the maximal time of existence. Then
\begin{eqnarray}
\sup_{x\in \R} u_x(t,x)\leq\sup_{x\in \R} u_{0,x}(x)+\big(2\sup_{x\in \R} \rho_{0}^2(x)+\frac{3}{8}+\frac{11\sqrt{6}}{36}\big)t,\quad\quad\quad\quad\quad\quad\label{5.3}\\
\inf_{x\in \R} u_x(t,x)\geq\inf_{x\in \R} u_{0,x}(x)+\big((\frac{1}{2}-\frac{3}{2}(E(0)+\varepsilon_0))\inf_{x\in \R} \rho_{0}^2(x)-\frac{1}{24\varepsilon_0}-\frac{7}{6}-\frac{11\sqrt{6}}{36}\big)t,\label{5.4}
\end{eqnarray}
where $0<\varepsilon_0<\frac{1}{3}-E(0).$
\end{lemma5}
\begin{proof}
Note that the local well-posedness result and a density argument indicates that it
suffices to prove the desired estimates for $s\geq 3$. We may assume that $u_0\not\equiv0$, otherwise the results are trivial. Since now $s\geq 3$, we have $u\in C^1_{0}(\R)$,
where the subscript $0$ means that the function decays to zero at infinity. Thus
\begin{eqnarray}\label{5.5}
\inf_{x\in \R} u_x(t,x)\leq0,\quad \mbox{and} \quad \sup_{x\in \R} u_x(t,x)\geq0, \quad t\in[0,T).
\end{eqnarray}

Differentiating the first equation of (\ref{two-componentnew}) with respect to $x$ and using the identity $(1-\partial_x^2)^{-1}f=p*f$, we get
\begin{eqnarray}\label{5.6}
{u_{tx}} =\frac{3}{2}u^2 +\frac{1}{2}\rho^2-3u^2(\rho^2-\rho)+4u^3-p \ast ( \frac{3}{2}u^2 +\frac{1}{2}\rho^2-3u^2(\rho^2-\rho)+4u^3).
\end{eqnarray}
Setting $M(t):=\sup\limits_{x\in \R}\{u_x(t,x)\},$ it then follows from Lemma $\ref{Lem5.1}$ that there exists $\gamma(t)\in \R$ such that $M(t)=u_x(t,\gamma(t)), t\in [0,T)$ and $u_{xx}(t,\gamma(t))=0$, a.e. $t\in [0,T)$.
Thanks to $q(t,\cdot)$ defined in (\ref{flow}) is a diffeomorphism of $\R$ for every $t\in[0,T)$, we deduce that there exists $x_1(t)\in \R$ such that $q(t,x_1(t))=\gamma(t), t\in [0,T).$ Now we define
$\bar{u}(t):=u(t,q(t,x_1(t)))$ and $\bar{\zeta}(t):=\rho(t,q(t,x_1(t))), t\in [0,T).$ Therefore, along the trajectory $q(t,x_1)$, we derive from (\ref{5.6}) and the second equation of (\ref{two-componentnew}) that
\begin{equation}\label{alongflow}
\begin{aligned}
\left\{ {\begin{array}{*{20}{l}}
\frac{d}{dt}M(t)=\frac{1}{2}\bar{\zeta}^2-3(\bar{\zeta}^2-\bar{\zeta})\bar{u}^2+f(t,q(t,x_1)),\\
\frac{d}{dt}\bar{\zeta}(t)=-\bar{\zeta}(t) M(t)(1-3\bar{u}^2),
\end{array}} \right.
\end{aligned}
\end{equation}
where $f:=\frac{3}{2}u^2+4u^3- p \ast ( \frac{3}{2}u^2 +\frac{1}{2}\rho^2-3u^2(\rho^2-\rho)+4u^3).$

Solving the second equation in (\ref{alongflow}), we obtain
\begin{eqnarray}\label{5.8}
\bar{\zeta}(t)=\bar{\zeta}(0)\exp\{-\int_0^t(M(1-3\bar{u}^2))(\tau)d\tau\}.
\end{eqnarray}
Taking account of the assumption on $E(0)$,
we get
\begin{eqnarray}\label{5.9}
\bar{u}^2(t)\leq \frac{1}{2}E(0)<\frac{1}{6}.
\end{eqnarray}
From (\ref{5.5}) we know that $M(t)\geq 0, t\in[0,T)$. Then combining (\ref{5.8})-(\ref{5.9}), we have
\begin{eqnarray}\label{5.10}
|\rho(t,q(t,x_1))|=|\bar{\zeta}(t)|\leq|\bar{\zeta}(0)|.
\end{eqnarray}

Now we estimate the upper bound of $f$. In view of the definition of $f$, we have
\begin{eqnarray*}
f \leq \frac{3}{2}u^2+4|u|^3+4|p \ast u^3|+3|p*u^2(\rho-1)^2|+3|p*u^2(\rho-1)|.
\end{eqnarray*}
Then it follows from the the conservation law, Young's inequality and H\"{o}lder inequality that
\begin{eqnarray*}
u^2\leq \frac{1}{2}\int_{\R}(u^2+u^2_x)dx\leq \frac{1}{2}E(0), \quad |u|^3\leq \frac{\sqrt{2}}{4} E^\frac{3}{2}(0),\\
|p \ast u^3|\leq  \|p\|_{L^\infty} \|u^3\|_{L^1}\leq \|p\|_{L^\infty}\|u\|_{L^\infty}\|u\|^2_{L^2}
\leq\frac{\sqrt{2}}{4} E^\frac{3}{2}(0),\\
|p*u^2(\rho-1)^2|\leq \|p\|_{L^\infty}\|u\|^2_{L^\infty}\|\rho-1\|^2_{L^2}\leq \frac{1}{4} E^2(0),
\end{eqnarray*}and
\begin{eqnarray*}
|p*u^2(\rho-1)|\leq \|p\|_{L^\infty} \|u\|_{L^\infty}\|u(\rho-1)\|_{L^1}\leq \frac{\sqrt{2}}{4} E^\frac{1}{2}(0)\|u\|_{L^2}\|\rho-1\|_{L^2}\leq \frac{\sqrt{2}}{4} E^\frac{3}{2}(0).
\end{eqnarray*}
Consequently, we obtain the upper bound of $f$
\begin{eqnarray}\label{upper}
f \leq\frac{3}{4}E(0)+\frac{11\sqrt{2}}{4} E^\frac{3}{2}(0)+\frac{3}{4} E^2(0).
\end{eqnarray}
Therefore, it follows from the first equation in (\ref{alongflow}), and (\ref{5.9})-(\ref{upper}) that
\begin{eqnarray*}
\frac{d}{dt}M(t)\leq 2 \bar{\zeta}^2 +\frac{3}{2}\bar{u}^4+f\leq 2\bar{\zeta}^2(0)+\frac{3}{8}+\frac{11\sqrt{6}}{36}.
\end{eqnarray*}
Integrating the above inequality from $0$ to $t$ yields the desired result (\ref{5.3}).

Now we turn to prove (\ref{5.4}). Setting $m(t):=\inf\limits_{x\in \R}\{u_x(t,x)\},$ there exists $\xi(t),x_2(t)\in \R$ such that $m(t)=u_x(t,\xi(t)),$ $q(t,x_2(t))=\xi(t), t\in [0,T),$
 and $u_{xx}(t,\xi(t))=0$, a.e. $t\in [0,T)$. Defining
$\tilde{u}(t):=u(t,q(t,x_2(t)))$ and $\tilde{\zeta}(t):=\rho(t,q(t,x_2(t))), t\in [0,T),$ along the trajectory $q(t,x_2)$, we get
\begin{equation*}
\begin{aligned}
\left\{ {\begin{array}{*{20}{l}}
\frac{d}{dt}m(t)=\frac{1}{2}\tilde{\zeta}^2-3(\tilde{\zeta}^2-\tilde{\zeta})\tilde{u}^2+f(t,q(t,x_2)),\\
\frac{d}{dt}\tilde{\zeta}(t)=-\tilde{\zeta}(t) m(t)(1-3\tilde{u}^2).
\end{array}} \right.
\end{aligned}
\end{equation*}
Since $m(t)\leq 0, t\in[0,T)$, a similar argument as (\ref{5.10}) implies
\begin{eqnarray}\label{5.12}
|\rho(t,q(t,x_2))|=|\tilde{\zeta}(t)|\geq|\tilde{\zeta}(0)|.
\end{eqnarray}
Similarly, we get the lower bound of $f$
\begin{eqnarray}\label{lower}
-f &\leq& 4|u|^3+\frac{3}{2}|p \ast u^2|+4|p \ast u^3|+\frac{1}{2}|p \ast (\rho-1)^2|+|p \ast (\rho-1)|\nonumber\\
&&+\frac{1}{2}p \ast1
+3|p*u^2(\rho-1)|\nonumber\\
&\leq&\frac{3}{4}+\frac{5}{4} E(0)+\frac{11\sqrt{2}}{4} E^\frac{3}{2}(0)\leq \frac{7}{6}+\frac{11\sqrt{6}}{36},
\end{eqnarray}
where we have used the inequalities $p*h^2\leq \frac{1}{2}\|h^2\|_{L^1}=\frac{1}{2}\|h\|^2_{L^2},$
and $|p \ast (\rho-1)|\leq\|p\|_{L^2}\|\rho-1\|_{L^2}=\frac{1}{2} \|\rho-1\|_{L^2}\leq \frac{1}{4} + \frac{1}{4}\|\rho-1\|^2_{L^2}.$
Then combining Young's inequality with (\ref{5.9}), (\ref{5.12})-(\ref{lower}), we obtain
\begin{eqnarray}\label{m}
\frac{d}{dt}m(t)&\geq&
(\frac{1}{2}- 3\tilde{u}^2-\frac{3\varepsilon_0}{2})\tilde{\zeta}(t)^2-\frac{3}{2\varepsilon_0}
\tilde{u}^4-\frac{7}{6}-\frac{11\sqrt{6}}{36}\nonumber\\
&\geq&
(\frac{1}{2}- \frac{3 }{2}(E(0)+\varepsilon_0))\tilde{\zeta}(t)^2-\frac{3}{8\varepsilon_0}
E^2(0)-\frac{7}{6}-\frac{11\sqrt{6}}{36}\nonumber\\
&\geq&(\frac{1}{2}- \frac{3 }{2}(E(0)+\varepsilon_0))\tilde{\zeta}(0)^2-\frac{1}{24\varepsilon_0}-\frac{7}{6}-\frac{11\sqrt{6}}{36},
\end{eqnarray}
where we take $\varepsilon_0$ satisfying $0<\varepsilon_0<\frac{1}{3}-E(0).$ Thus integrating (\ref{m}) over $[0,t]$ implies (\ref{5.4}).
This completes the proof of Lemma $\ref{Lem5.2}$.
\end{proof}
Finally, we are in the position to establish the global existence result for $(\ref{two-componentnew})$.
\begin{theorem5}\label{Th5.1}
Assume $E(0)<\frac{1}{3}$. If $(u_0,\rho_0-1)\in H^s\times H^{s-1}$, $s>\frac{3}{2}$, then there exists a unique solution $(u,\rho)$ of $(\ref{two-componentnew})$ in $C([0,\infty);H^s\times H^{s-1})\cap C^1([0,\infty);H^{s-1}\times H^{s-2})$
with the initial data $(u_0,\rho_0)$. Moreover, the solution $(u,\rho)$ depends continuously on initial data $(u_0,\rho_0)$ and the conserved quantity $E(t)$ is independent of the existence time.
\end{theorem5}
\begin{proof}
With the aid of Lemma $\ref{Lem5.2}$, we thus deduce from the blow-up criterion in Remark $\ref{Rem4.1}$ that the local solution $(u,\rho)$ guaranteed by Remark $\ref{Rem3.1}$ can be extended to all of the interval $[0,\infty).$
This completes the proof of Theorem \ref{Th5.1}.
\end{proof}

\section{Appendix}
\newtheorem {remark6}{Remark}[section]
\newtheorem{theorem6}{Theorem}[section]
\newtheorem{definition6}{Definition}[section]
\newtheorem{proposition6}{Proposition}[section]
\newtheorem{lemma6}{Lemma}[section]
In Appendix, we recall some basic theory of the
Littlewood-Paley decomposition and Besov spaces for completeness. One can check \cite{Bahouri,Danchin} for details.
There exist two smooth radial functions $\chi(\xi)$ and $\varphi(\xi)$ valued in $[0,1]$, such that $\chi$ is supported in
$\mathcal{B}=\{\xi\in\R^d,|\xi|\leq\frac{4}{3}\}$ and $\varphi(\xi)$ is supported in $\mathcal{C}=
\{\xi\in\R^d,\frac{3}{4}\leq|\xi|\leq\frac{8}{3}\}$. Denote $\mathcal{F}$ by the Fourier transform and $\mathcal{F}^{-1}$ by its inverse.
For all $u\in \mathcal{S'}(\R^d)$($\mathcal{S'}(\R^d)$ denotes the tempered distribution spaces), the nonhomogeneous dyadic operators $\Delta_q$ and
the low frequency cut-off operator $S_q$ are defined as follows:
$ \Delta_q u=0$ for $q\leq-2$, $\Delta_{-1}u=\chi(D)u=\mathcal{F}^{-1}(\chi\mathcal{F}u),$
$\Delta_qu=\varphi(2^{-q}D)u=\mathcal{F}^{-1}(\varphi(2^{-q}\cdot)\mathcal{F}u),
$ for $q\geq0,$ and $S_q u=\sum^{q-1}_{i=-1}\Delta_i u=\chi(2^{-q}D)u=\mathcal{F}^{-1}(\chi(2^{-q}\cdot)\mathcal{F}u).$
\begin{definition6}\label{Def6.1}
 (Besov spaces) Assume $s\in \R,\ 1\leq p,r\leq\infty$. The nonhomogeneous Besov
 space $B^s_{p,r}(\mathbb{R}^d)$ ($B^s_{p,r}$ for short) is defined by
\begin{equation*} B^s_{p,r}=\{u\in\mathcal{S'}(\mathbb{R}^d): \|u\|_{B^s_{p,r}}=\|2^{qs}\Delta_qu\|_{l^r(L^p)}=\|(2^{qs}\|\Delta_qu\|_{L^p})_{q\geq -1}\|_{l^r}<\infty\}.\end{equation*}
In particular, $B^\infty_{p,r}=\bigcap_{s\in
\R}B^s_{p,r}.$
\end{definition6}

Now we list some useful properties and lemmas about Besov space which are used in the previous sections.
 \begin{lemma6}\label{Lem6.1}
Let $s\in \R,\ 1\leq p,r,p_i,r_i\leq\infty,i=1,2$. Then\\
(i) Algebraic properties: if $s>0$, $B^s_{p,r}\cap L^\infty$ is a Banach
algebra. Moreover, $B^s_{p,r}$ is an algebra $\Leftrightarrow B^s_{p,r}\hookrightarrow L^\infty$
$\Leftrightarrow s>\frac{d}{p}$ or $s\geq\frac{d}{p}$ and $r=1.$\\
(ii) Fatou lemma: if $\{u_n\}_{n\in \mathbb{N}}$ is a bounded sequence in
$B^s_{p,r}$, then there exist an element $u \in B^s_{p,r}$ and a subsequence $\{u_{n_k}\}_{k\in \mathbb{N}}$ such that\\
$$\lim_{k\rightarrow \infty}u_{n_k}= u \ \mbox{in} \ \mathcal{S'}, \quad \mbox{and} \quad \|u\|_{B^s_{p,r}}\leq \liminf_{k\rightarrow \infty}\|u_{n_k}\|_{B^s_{p,r}}.$$
(iii) Interpolation: (1) if $u\in B^{s_1}_{p,r}\cap B^{s_2}_{p,r}$ and $\theta\in[0,1]$, then $u\in B^{\theta s_1+(1-\theta)s_2}_{p,r}$ and
$\|u\|_{B^{\theta s_1+(1-\theta)s_2}_{p,r}}\leq
\|u\|^\theta_{B^{s_1}_{p,r}}\|u\|^{1-\theta}_{B^{s_2}_{p,r}}.$
(2) if $u\in B^{s_1}_{p,\infty} \cap B^{s_2}_{p,\infty}$ and $s_1<s_2$, then $u\in B^{\theta s_1+(1-\theta)s_2}_{p,1}$ for all
$\theta\in(0,1)$ and there exists a constant $C$ such that
$\|u\|_{B^{\theta s_1+(1-\theta)s_2}_{p,1}}\leq \frac{C}{\theta(1-\theta)(s_2-s_1)}
\|u\|^\theta_{B^{s_1}_{p,\infty}}\|u\|^{1-\theta}_{B^{s_2}_{p,\infty}}.$\\
(iv) Action of Fourier multipliers on Besov spaces: let $m\in \R$
and $f$ be a $S^m$-multiplier ($i.e.,$ $f:\R^d\rightarrow \R$ is a
smooth function and satisfies that for each multi-index $\alpha$,
there exists a constant $C_\alpha$ such that $|\partial^\alpha
f(\xi)|\leq C_\alpha(1+|\xi|)^{m-|\alpha|}$, for $\forall \xi \in
\R^d$). Then the operator $f(D)=\mathcal{F}^{-1}(f\mathcal{F})\in Op(S^{m})$ is continuous from $B^s_{p,r}$ to
$B^{s-m}_{p,r}$.
\end{lemma6}

\begin{lemma6}\label{Lem6.2} (A priori estimate)
Let $1\leq p,r\leq +\infty, s>-\min\{\frac{1}{p},1-\frac{1}{p}\}$. Assume that $f_0\in B^s_{p,r}$, $F\in L^1(0,T;B^s_{p,r})$, and
$\partial_x v$ belongs to $L^1(0,T;B^{s-1}_{p,r})$ if $s>1+\frac{1}{p}$
 or to $L^1(0,T;B^{\frac{1}{p}}_{p,r}\cap L^\infty)$ otherwise. If $f\in L^\infty(0,T;B^{s}_{p,r})\cap
C([0,T];\mathcal{S'})$ solves the following transport equation
 \begin{equation}\label{Eq5.1}
\left\{\begin{array}{ll}\partial_tf+v\partial_x f=F,\\
 f\big|_{t=0}=f_0,\\
\end{array}\right.
\end{equation}
then there exists a constant $C$ depending only on $s,p,r$ such that
the following statements hold for $t\in [0,T]$\\
(i) If $r=1$ or $s\neq 1+\frac{1}{p}$,
\begin{eqnarray*}
\|f(t)\|_{B^{s}_{p,r}}\leq
\|f_0\|_{B^{s}_{p,r}}+\int_0^t\|F(\tau)\|_{B^{s}_{p,r}}d\tau
+C\int_0^tV'(\tau)\|f(\tau)\|_{B^{s}_{p,r}}d\tau,
\end{eqnarray*}
 or
 \begin{eqnarray*}
\|f(t)\|_{B^{s}_{p,r}}\leq
e^{CV(t)}(\|f_0\|_{B^{s}_{p,r}}+\int_0^te^{-CV(\tau)}
\|F(\tau)\|_{B^{s}_{p,r}}d\tau),
\end{eqnarray*}
where $V(t)=\int_0^t\|\partial_xv(\tau,\cdot)\|_{B^{\frac{1}{p}}_{p,r}\cap L^\infty}d\tau$ if $s<1+\frac{1}{p}$, and $V(t)=\int_0^t\|\partial_xv(\tau,\cdot)\|_{B^{s-1}_{p,r}}d\tau$ otherwise.\\
(ii) If $r<\infty$, then $f\in C([0,T];B^{s}_{p,r})$. If $r=\infty$,
then $f\in C([0,T];B^{s'}_{p,1})$ for all $s'<s$.
\end{lemma6}

\begin{lemma6}\label{Lem6.3}
(Existence and uniqueness) Let $p,r,s,f_0$ and $F$ be as in the
statement of Lemma \ref{Lem6.2}. Suppose that $v\in
L^\rho(0,T;B_{\infty,\infty}^{-M})$ for some $\rho>1,M>0$ and
$\partial_x v\in L^1(0,T;B^{\frac{1}{p}}_{p,\infty}\cap L^\infty)$
if $s<1+\frac{1}{p}$, and $\partial_x v\in L^1(0,T;B^{s-1}_{p,r})$
if $s>1+\frac{1}{p}$ or $s=1+\frac{1}{p}$ and $r=1$. Then (\ref{Eq5.1}) has a unique solution $f\in
L^\infty(0,T;B^{s}_{p,r})\cap
(\cap_{s'<s}C([0,T];B^{s'}_{p,1}))$ and the corresponding
inequalities in Lemma \ref{Lem6.2} hold true. Moreover, if
$r<\infty$, then $f\in C([0,T];B^{s}_{p,r}).$
\end{lemma6}

\begin{lemma6}\label{Lem6.4}
 (Commutator estimates) Let $s>0$, $1\leq p,r \leq \infty$, and $v$ be a vector field over $\R^d$. Define $R_q=[v\cdot\nabla,\Delta_q]f$, where $[\cdot,\cdot]$ denotes the commutator of the operators. Then there exists a constant $C$
 such that
 \begin{eqnarray*}
\|(2^{qs}\|R_q\|_{L^p})_{q\geq -1}\|_{l^r}\leq C (\|\nabla v\|_{L^\infty}\|f\|_{B_{p,r}^s}+\|\nabla f\|_{L^\infty}\|\nabla v\|_{B_{p,r}^{s-1}}).
\end{eqnarray*}
\end{lemma6}

\noindent\textbf{Acknowledgments} \
The work is supported by National Nature Science Foundation of China under Grant 12001528. The authors thank the anonymous
referee for helpful suggestions and comments.

\end{document}